\documentclass[graybox,envcountsame,envcountsect,numbook]{svmult}
\usepackage{mathptmx}       
\usepackage{helvet}         
\usepackage{courier}        
\usepackage{makeidx}         
\usepackage{float}           
\usepackage{graphicx}        
\usepackage{multicol}        
\usepackage[bottom]{footmisc}
\usepackage{amsmath,amssymb,amsfonts}        
\makeindex             

\usepackage{amsmath, latexsym, amsfonts,bm,amssymb,mathscinet} 
\usepackage{mathtools}
\usepackage{xcolor}
\usepackage[ruled,vlined]{algorithm2e}
 
\usepackage{url} 
\usepackage{enumerate}
\usepackage{units} 

\newtheorem{cnd}{Condition}

\newtheorem{assump}[theorem]{Assumption}

\def\exzb{{ {\mathbb{E}_{a_0,s_0}\,}}}

\newcommand{\ind}{\mathbf{1}}

\def\ee{{ {\mathbb{E}}}}
\def\eps{\varepsilon}

\renewcommand{\th}{\theta}

\newcommand{\ig}{\operatorname{IG}}

\def\epsilon{\varepsilon}

\newcommand{\cH}{{\mathcal{H}}}

\newcommand{\dd}{\mathrm{d}}
\newcommand{\pp}{\mathbb{P}}
\newcommand{\fpspace}{(\Omega,\mathcal{F},\mathbb{F},\pp)}
\newcommand{\cf}{\mathcal{F}}
\newcommand{\one}{\mathbf{1}}
\renewcommand\thesection{\thechapter.\arabic{section}}

\newcommand{\EE}{\mathbb{E}}
\newcommand{\VV}{\mathbb{V}}

\def\iid{\stackrel{\textrm{i.i.d.}}{\sim}}

\newcommand{\g}{\operatorname{G}}

\allowdisplaybreaks

\begin{document}

\title*{Nonparametric {B}ayesian inference for stochastic processes with piecewise constant priors}
\author{Denis Belomestny, Frank van der Meulen, and Peter Spreij}
\authorrunning{Belomestny, van der Meulen, Spreij}
\institute{Denis Belomestny \at 
Faculty of Mathematics,
Duisburg-Essen University,
Thea-Leymann-Str.~9,
D-45127 Essen,
Germany
and Faculty of Computer Sciences,
HSE University,
Moscow, Russian Federation.
\email{denis.belomestny@uni-due.de}
\and Frank van der Meulen \at Department of Mathematics,
Vrije Universiteit,
De Boelelaan 1111,  
1081HV Amsterdam,
The Netherlands. 
\email{f.h.van.der.meulen@vu.nl}
\and Peter Spreij \at Korteweg-de Vries Institute for Mathematics,
University of Amsterdam,
P.O. Box 94248,
1090 GE Amsterdam,
The Netherlands and Institute for Mathematics, Astrophysics and Particle Physics, Radboud University, Nijmegen, The Netherlands.~\email{spreij@uva.nl}}

\titlerunning{Nonparametric {B}ayesian inference for stochastic  with piecewise constant priors}
\maketitle

\abstract{We present a survey of some of our recent results on Bayesian nonparametric inference for a multitude of stochastic processes. The common feature is that the prior distribution in the cases considered is on suitable sets of piecewise constant or piecewise linear functions, that differ for the specific situations at hand. Posterior consistency and in most cases contraction rates for the estimators are presented. Numerical studies on simulated and real data accompany the theoretical results.
}


\section{Introduction}

This paper is a survey of our recent results, obtained with collaborators, in Bayesian nonparametric function estimation. Although the settings and corresponding models that we consider exhibit a large variety (volatility estimation
in diffusion models, Poisson intensity estimation, estimation for gamma-driven stochastic differential equations, estimation in the presence of microstructure noise, and L\'evy density estimation for Gamma-type subordinators), they all share the following common property. The function to be estimated, which can be a function of time or  space, is approximated by piecewise constant functions (with one exception that deals with piecewise linear functions). This assumption is one of the possibilities to reduce the estimation of an infinite dimensional parameter to the estimation of only finitely many parameters, this being inevitable if one has to \emph{compute} estimators. 

We adopt the Bayesian approach by specifying a prior distribution supported on piecewise constant functions. The goal is then to compute the posterior distribution and show it has good properties, when considered from a frequentist point of view. In this work, the first goal is either accomplished by using conjugate priors or using Markov Chain Monte Carlo to obtain samples from the posterior. The second goal is achieved by showing posterior consistency and quantification of posterior contraction rates. A thorough introduction to these concepts is given in \cite{ghosal17}. Informally, posterior consistency means that the posterior mass of a fixed $\epsilon$-ball around the ``true'' function tends to $1$ as the sample size increases. The contraction rates quantifies the rate at which we can shrink the radius of these balls with increasing sample size while still asymptotically having posterior mass $1$.  

It turns out that this procedure is feasible in all cases that we consider. It is intuitively clear that rougher (continuous) functions are more difficult to estimate than smoother ones. The degree of smoothness is measured by the H\"older exponent in our approach. We assume in all cases that the functions to be estimated are H\"older continuous. We obtain contraction rates that depend on the H\"older exponent, and the smoother the function, the better the estimation procedure's contraction rate.
As we will demonstrate, despite the entirely different settings, the contraction rates are identical in all cases. 

For theoretical and numerical implementations, one has to specify priors on the approximating piecewise continuous functions. We will assume a fixed set of bins, reflecting the resolution at which we aim to learn the unknown function. In a first approach, we will a priori take the (finitely many) values on the bins as \emph{independent} random variables with  distributions from an appropriate chosen class (mostly gamma or inverse gamma distributions), such that the prior is (partially) conjugate. These specifications lead to the theoretical results on the posterior distributions mentioned above. 

It also appeared that better numerical results, clearly visible in the figures we provide and less sensitive to the chosen set of bins, can be obtained by relaxing the independence assumption on the prior distribution. We accomplish this observation by using  appropriate Markov chain priors on the levels of the piecewise constant functions. A principal motivation is that such a prior introduces dependence between the adjacent levels of the piecewise constant approximating functions, which is a rather natural property if one realizes that values of continuous functions tend to be close to each other on adjacent intervals.
The theoretical results on the contraction rates we show in several cases coincide with the minimax convergence rates in frequentist estimation. From that asymptotic point of view, the estimators with the Markov chain prior cannot be better than the ones with the independent priors. Nevertheless, in practical situations, they turn out to be more accurate.

A great variety of numerical examples in all settings are presented. They include results from simulated examples to illustrate the theoretical properties and from real data to derive some conclusions based on our methods.

The paper is organized as follows. Section~\ref{section:diffusion} details the results on volatility estimation
for stochastic differential equations driven by a Brownian motion, Section~\ref{section:micro} is devoted to estimation for partial observations of diffusion processes in presence of microstructure noise,
Section~\ref{section:poisson} treats intensity estimation for Poisson processes,
Section~\ref{section:gammasde} gives results on estimation of the volatility function in gamma-driven stochastic differential equations,
Section~\ref{section:gamma} is on estimation of the L\'evy density for Gamma-type  subordinators.
Finally Section~\ref{section:conclusion} summarizes and concludes. 

\medskip

A remark on the notation throughout the paper: for two sequences $(a_n)$ and $(b_n)$ of real numbers we write $a_n\asymp b_n$ if there are constants $c$ and $C$ in $(0,\infty)$ such that eventually $cb_n\leq a_n\leq Cb_n$.

\section{Volatility estimation}\label{section:diffusion}

In this section, we consider the nonparametric Bayesian estimation of a volatility function, also synonymously referred to as the dispersion coefficient, that arises in diffusion models. This function may be a smooth function of the time parameter, but it may also happen that it is a realization of another diffusion process, like a Brownian motion. Hence the function to be estimated may have various orders of regularity, which we characterize by their degree of H\"older continuity. 
This section is based on \cite{gugu17}, \cite{gugu18b}, which also contain many references to related literature.

\subsection{Problem formulation}

Stochastic differential equations (SDEs) driven by a Brownian motion have been widely used as models in numerous applications ranging from physics to engineering and to finance. Their solutions, weak or strong, are diffusion processes. In this section we consider a one-dimensional diffusion process  $X$  satisfying the stochastic differential equation (SDE)
\begin{equation}
	\label{sde}
	\dd X_t=a(t,X_t)\,\dd t+s(t)\,\dd W_t, \quad X_0=x_0, \quad t\in[0,1].
\end{equation}
Here $W$ is a standard Wiener process, and $a$ and $s$ are the  drift function and volatility function, respectively. Note that $s$ is a function of time only. It is further assumed that $b$ and $s$ are such that \eqref{sde} admits a weak solution that is unique in law.
We assume that observations 
\[\scr{X}_n =\{X_{t_{0,n}}, \ldots, X_{t_{n,n}}\}\] 
from the solution $X$ to \eqref{sde} are available, where $t_{i,n}=\tfrac{i}n,$ $i=0,\ldots,n.$ In the statistical set-up, we denote the true drift and diffusion coefficient by $a_0$ and $s_0$,
Our aim is to estimate $s_0$ nonparametrically within the Bayesian setup. We note that nonparametric Bayesian estimation of volatility is very different from the drift coefficient estimation: the latter fundamentally relies on the equivalence of laws of continuously observed diffusion processes that have the same diffusion coefficient. This is not applicable when the diffusion coefficient itself is unknown and is a parameter to be estimated when  a continuous record of observations is available. In the current set-up on the contrary, with \emph{discrete time} observations, there is  equivalence of laws and it is therefore possible to derive a likelihood, at least theoretically, a fundamental tool for Bayesian inference. An implementable version of such a likelihood is obtained by \emph{intentionally} misspecifying the drift coefficient $a$ by taking it equal to zero, and employing a (conjugate) histogram-type prior on the diffusion coefficient, that has piecewise constant realisations on bins forming a partition of $[0,1]$. Due to this, our nonparametric Bayesian method to estimate the dispersion coefficient $s_0$ in \eqref{sde}  is  easily  implemented. Theoretical asymptotic results for estimating $s$ will be obtained in the case when $n\to\infty$ (high frequency observations), keeping the observation interval $[0,1]$ fixed. For simplicity, we take  $x_0=0$. Because of the fixed interval $[0,1]$, consistent estimation of the drift is impossible, but the high frequency observations make it possible to ignore the drift if the aim is to (consistently) estimate the volatility. We will make this precise in Proposition~\ref{prop:ignore}.

\subsection{Likelihood, prior and posterior}\label{section:lpp}

Under the assumption that the drift in \eqref{sde} is zero, the (pseudo-)likelihood is that of a Gaussian sample and is given by
\begin{equation}
\label{likelihvol}
L_n(s)=\prod_{i=1}^{n} \left\{ \frac{1}{\sqrt{2\pi\int_{t_{i-1,n}}^{t_{i,n}}s^2(u)\,\dd u}}\phi_0\left( \frac{Y_{i,n}}{\sqrt{\int_{t_{i-1,n}}^{t_{i,n}}s^2(u)\,\dd u}} \right) \right\},
\end{equation}
where $Y_{i,n}=X_{t_{i,n}}-X_{t_{i-1,n}}$ and
$\phi_0(u)=\exp(-u^2/2).$ With $\Pi_n$ denoting a prior on the dispersion coefficient, provided all the involved quantities are suitably measurable, Bayes' theorem gives that the posterior probability of any measurable subset $S\subset\mathcal{S}$ of dispersion coefficients is given by
\begin{equation*}
\Pi_n(S\mid X_{t_{0,n}}\ldots,X_{t_{n,n}})=\frac{\int_{S} L_n(s) \Pi_n(\dd s)}{ \int_{\mathcal{S}} L_n(s) \Pi_n(\dd s) },
\end{equation*}
where $\mathcal{S}$ denotes a space on which the prior $\Pi_n$ is defined. 

It follows from \eqref{likelihvol} that the likelihood depends on the parameter of interest only through the integrals $\int_{t_{i-1,n}}^{t_{i,n}}s^2(u)\,\dd u$. Consequently, it appears natural to a priori model the diffusion coefficient as piecewise constant on intervals $[{t_{i-1,n}},{t_{i,n}}]$. However, some smoothing should also be performed, and this can be achieved by aggregated several neighbouring intervals $[{t_{i-1,n}},{t_{i,n}}]$. Thus, to construct the prior, we proceed as follows: Let $m$ be an integer smaller than $n$. 
Then we can uniquely write $n=mN+r$  with $0\leq r<m$, and in fact $N=\lfloor \frac{n}{m}\rfloor$. 
 Both $m$ and $N$ will depend on $n$ (and we also write $m_n$ and $N_n$ to emphasize this when appropriate). With this assumption we have bins $B_k=[t_{m(k-1),n},t_{mk,n})$, $k=1,\ldots,N-1$ and $B_N=[t_{m(N-1),n},1]$. 
 
Let $\mathcal{S}_n$ be the set of functions $s$ that can be represented as 
$
s=\sum_{k=1}^{N_n} \xi_k \ind_{B_k}.
$
The prior $\Pi_n$ on the dispersion coefficient $s\in\mathcal{S}_n$ is defined by putting a prior on the coefficients $\xi_k$'s. Since
\begin{equation}
\label{eq:series}
s^2=\sum_{k=1}^{N_n} \xi_k^2 \ind_{B_k}=\sum_{k=1}^{N_n} \theta_k \ind_{B_k},
\end{equation}
where we have put $\theta_k=\xi_k^2$, equivalently one can place the prior on the coefficients $\theta_k$'s of the diffusion coefficient $s^2$. We call the prior $\Pi_n$ a histogram-type prior. 

As a result of choosing $\Pi_n$ such that the $\theta_i$ are independent having an inverse gamma distribution, we can explicitly compute the conjugate posterior of the $\theta_i$, again of inverse gamma type.
\begin{lemma}\label{lem:post-ig}
Assume $\th_1,\ldots, \th_N$ are independent with the inverse gamma $\mathrm{IG}(\alpha,\beta)$ distribution. Then $\th_1,\ldots, \th_N$ are a posteriori independent and, for $k=1,\ldots,N-1$,
\[ \th_k \mid \scr{X}_n \sim \mathrm{IG}\left(\alpha+m/2, \beta+ n Z_k/2\right). \]
with 
\[ Z_k = \sum_{i=(k-1)m+1}^{km} Y_{i,n}^2, \]
whereas
\[ \th_N \mid \scr{X}_n \sim \mathrm{IG}\left(\alpha+(m+r)/2, \beta+ n Z_N/2\right), \]
with 
\[ Z_N = \sum_{i=(N-1)m+1}^{n} Y_{i,n}^2. \]
\end{lemma}

\subsection{Asymptotics}\label{section:asympvol}
We let  $\pp_{a_0,s_0}$ denote the law of the path $(X_t\colon t\in[0,1])$ from \eqref{sde} under the true parameter values $(a_0,s_0)$. In particular the notation $\mathbb{P}_{0,s_0}$ is used for such a law when the drift coefficient is equal to zero. For expectation we use the same subscripts.
We denote the prior distribution on the dispersion coefficient by $\Pi_n$ (with $n$ the number of observations) and write the posterior as $\Pi_n(\,\cdot\mid \mathcal{X}_n)$. 
By $\|\cdot\|_2$ we denote  the $L_2$-norm with respect to the Lebesgue measure on the Borel sets of $[0,1]$. 

The standing assumption  on the volatility in this section is the following.
\begin{assump}
\label{standing}
\emph{
Assume that the true volatility function $s_0$ is H\"older continuous on $[0,1]$ with H\"older constant $L$ and H\"older exponent $\lambda\in (0,1]$, $|s_0(u)-s_0(v)|\leq L|u-v|^\lambda$ for all $u,v\in [0,1]$,
and is bounded away from zero.
}
\end{assump}
The first result formalises that under our asymptotic regime, we can safely ignore the drift by, erroneously in general, putting it equal to zero.

\begin{proposition}\label{prop:ignore}
\label{prop1}
Let, for $\varepsilon>0,$ 
\begin{equation*}
U_{s_0,\varepsilon}=\left\{ s\in\mathcal{S}_n\colon\|s-s_0\|_2<\varepsilon \right\}
\end{equation*}
be the  $L_2$-neighbourhood of $s_0$ of radius $\varepsilon$. 
Let Assumption~\ref{standing} hold and assume that there exists a sequence $\varepsilon_n\rightarrow 0$
\begin{equation*}
\ee_{0,s_0}[\Pi_n(U_{s_0,\varepsilon_n}^c\mid\mathcal{X}_n) ] \rightarrow  0
\end{equation*}
as $n\rightarrow\infty.$
Then also
\begin{equation*}
\ee_{a_0,s_0}[\Pi_n(U_{s_0,\varepsilon_n}^c\mid\mathcal{X}_n)] \rightarrow 0.
\end{equation*}
\end{proposition}
As we consider the asymptotics $n\rightarrow\infty$, we take the number of bins $N$ to depend on the sample size $n$, and indicate this in our notation by writing $N_n$. 
The next result gives the posterior contraction rate $n^{-\lambda/(2\lambda+1)}$ for the $L_2$-norm. Establishing posterior contraction in the $L_2$-metric is rather natural, as $\int_0^1 s_0^2(t)\dd t$ is the quadratic variation of the process $(X_t:t\in [0,1])$ over the interval $[0,1].$  

\begin{theorem}\label{thm:rate-posterior-ig} 
Assume that Assumption~\ref{standing} holds.  Assume  $N_n \asymp n^{1/(2\lambda+1)}$.  If we let  $ \eps_n\asymp n^{-\lambda/(2\lambda+1)}$, then for any sequence $h_n$ tending to infinity (as $n\to \infty$) we have
\[ 
\exzb \left[ \Pi_n(\|s^2-s_0^2\|_2 \ge h_n \eps_n \mid \scr{X}_n) \right]  \to 0
\]
as $n\to \infty.$
\end{theorem}
A comparison with the frequentist minimax convergence rate in \cite{hoffmann97} shows that the posterior for the diffusion coefficient contracts at the optimal rate in the $L_2$-metric. With histogram-type priors considered in this work no further improvement in the posterior contraction rate, when $\lambda=1$, is possible beyond $n^{-1/3},$ even if the function $s_0$ is smoother than a Lipschitz function. An intuitive reason for this is that realisations of our histogram-type priors are too rough for this; cf.\ p.~629 in \cite{scricciolo07}.

\subsection{Numerical results}

In this section we illustrate the theoretical results of Section~\ref{section:asympvol} with a number of numerical examples.
For an implementation of our Bayesian analysis a judicious choice of the number of bins $N$ is crucial. 
For obtaining such a choice there are various methodologies. For instance one can use the Deviance Information Criterion (DIC) of \cite{spiegelhalter02} to select $N$, or choose $N$ as the maximizer of the marginal likelihood
$
n\mapsto \int_{\mathcal{S}} L_n(s) \Pi_n(\dd s).
$
We refer to \cite{gugu17} for details.
We use simulations of diffusion processes with known drift and diffusion coefficient to gain insight into the numerical performance of our method. We are particularly interested in the practical consequence of using a pseudo-likelihood ignoring the drift. 

First we simulate realisations from the model for different dispersion and drift coefficients. Given subsamples of those realisations sampled at different rates we
compute the posterior distribution of the dispersion function using the piecewise constant prior (histogram-type prior) of Lemma~\ref{lem:post-ig} with varying numbers of bins. As an illustration, plots of  marginal posterior bands are compared with the true dispersion function.
The marginal posterior bands are obtained by computing $1-\alpha$ central posterior intervals (see \cite{gelman}, p.~33) separately for the coefficients $\theta_k$'s using Lemma~\ref{lem:post-ig}. This way of computing the posterior intervals will also be used in later sections.
  By Proposition~\ref{prop1}, assuming that there is no drift still leads to consistent Bayesian estimation of the dispersion coefficient,  even if the data are from a diffusion process with nonzero drift. This is illustrated by our simulation results of this section.

We simulated sample paths of the diffusion $(X_t:t \in [0,1])$ where the true dispersion coefficient is given by one of
\begin{align*}
s_1(t) &= 3/2 + \sin(2(4t-2)) + 2\exp(-16(4t-2)^2),\\
s_2(t) &= B_t(\omega_0) + 1,
\intertext{
and the drift is given by one of  }
a_1(x) &= 0, \\
a_2(x) &= -10x+20.
\end{align*}
Note the change of notation in these examples, true dispersion and drift are not denoted $s_0$ and $a_0$ here, but differently to distinguish these special choices from the general case.
The function $s_1$ is a benchmark function used in \cite{fan2018} in the context of nonparametric regression, up to a vertical shift to ensure positivity. To define $s_2(t)$, we took a fixed realisation of a Brownian motion $B$, indicated by $\omega_0$, starting at $1$  with $B_t(\omega_0)  > -1$ for $t \in [0,1]$. 
The function $s_1$ is Lipschitz continuous, while  $s_2$ is H\"older continuous with exponent essentially $\frac12$. As the prior on the coefficients on the individual bins we used  $\operatorname{IG(0.1, 0.1)}$ distributions and the coefficients were assumed to be independent.

Figures~\ref{fig:s2} and \ref{fig:s4}  show marginal $98\,\%$ posterior bands for different combinations of bin number and observation regime for both dispersion coefficients when the drift is zero. Figure~\ref{fig:s2drift} shows the marginal posterior bands for $s_1$ which are obtained if an affine drift term $a_2(x) = -10x+20$ is present, but neglected in the estimation procedure. Comparison with Figure~\ref{fig:s2} shows that presence of a strong nonzero drift
hardly affects the obtained credible bands. Note that credible bands successfully recover the overall shape of the functions $s_i$, although the recovery is not too refined. 
The recovery is somewhat less accurate for function $s_2$, see Figure~\ref{fig:s4}, than for function $s_1,$ see Figure~\ref{fig:s2}. This is in perfect agreement with our theoretical results from Theorem~\ref{thm:rate-posterior-ig}, that give a slower posterior contraction rate for less smooth functions.

\newcommand{\pd}{\hspace{0.3cm}} 
\renewcommand{\arraystretch}{0}
\begin{figure}[H]
\rotatebox{00}{
\begin{tabular}{cccc}
& $N= 40$ & $N = 80$ & $N = 160$\\
\rotatebox[origin=rt]{-90}{$n=8001$\pd}&
\includegraphics[height=0.19\textwidth]{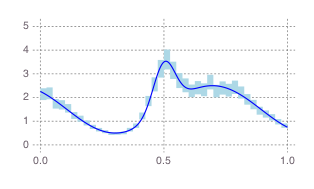} & 
\includegraphics[height=0.19\textwidth]{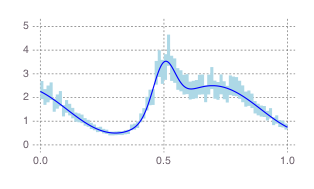} &
\includegraphics[height=0.19\textwidth]{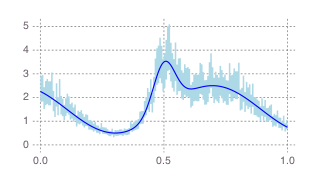} \\
\rotatebox[origin=rt]{-90}{$n=16\,001$\pd}&
\includegraphics[height=0.19\textwidth]{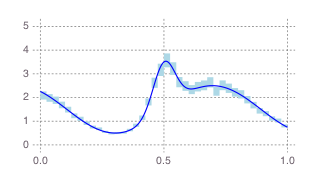}&
\includegraphics[height=0.19\textwidth]{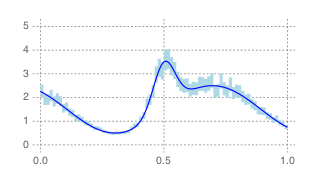}&
\includegraphics[height=0.19\textwidth]{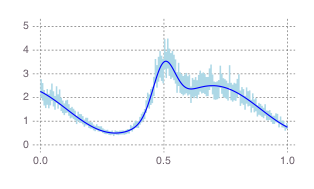}\\
\rotatebox[origin=rt]{-90}{$n=32\,001$\pd}&
\includegraphics[height=0.19\textwidth]{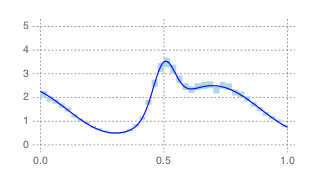}&
\includegraphics[height=0.19\textwidth]{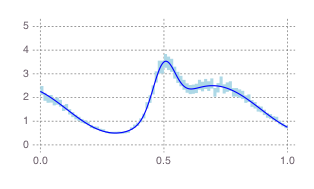}&
\includegraphics[height=0.19\textwidth]{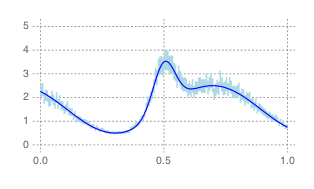}
\end{tabular}
}
\caption{Estimation results for $s_1$ with varying number of observations and bins, no drift.
The solid curve is the true dispersion coefficient while the light blue areas are $98\,\%$ marginal posterior bands. $n$ is the number of observations, $N$ is the number of bins.
}
\label{fig:s2}
\end{figure}	

\begin{figure}
\rotatebox{00}{
\begin{tabular}{cccc}
& $N= 40$ & $N = 80$ & $N = 160$\\
\rotatebox[origin=rt]{-90}{$n=8001$\pd}&
\includegraphics[height=0.19\textwidth]{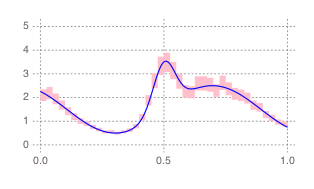} & 
\includegraphics[height=0.19\textwidth]{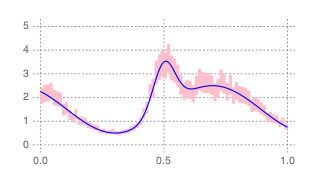} &
\includegraphics[height=0.19\textwidth]{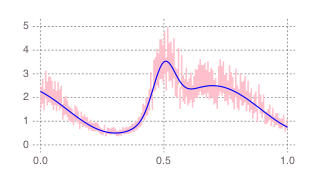} \\
\rotatebox[origin=rt]{-90}{$n=16\,001$\pd}&
\includegraphics[height=0.19\textwidth]{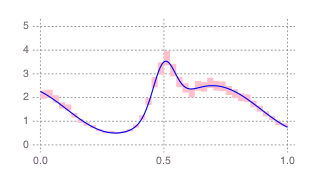}&
\includegraphics[height=0.19\textwidth]{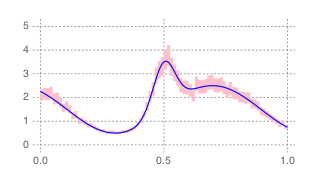}&
\includegraphics[height=0.19\textwidth]{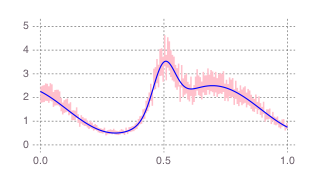}\\
\rotatebox[origin=rt]{-90}{$n=32\,001$\pd}&
\includegraphics[height=0.19\textwidth]{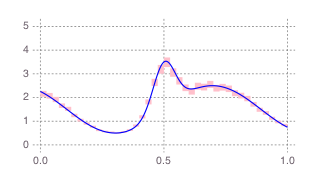}&
\includegraphics[height=0.19\textwidth]{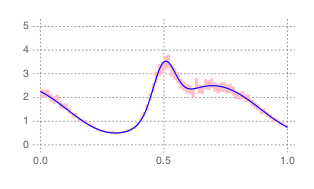}&
\includegraphics[height=0.19\textwidth]{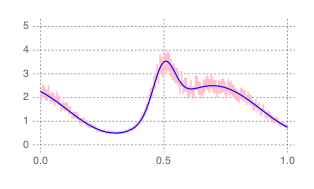}
\end{tabular}
}
\caption{Estimation for $s_1$ with varying number of observations and bins and drift $a_2(x) =-10x + 20$.
The solid curve is the true dispersion coefficient while the light red areas are $98\,\%$ marginal posterior bands. $n$ is the number of observations, $N$ is the number of bins.
}
\label{fig:s2drift}
\end{figure}

\begin{figure}
\rotatebox{00}{
\begin{tabular}{cccc}
& $N= 40$ & $N = 80$ & $N = 160$\\
\rotatebox[origin=rt]{-90}{$n=8001$\pd}&
\includegraphics[height=0.19\textwidth]{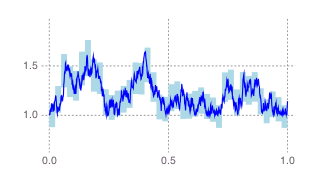} & 
\includegraphics[height=0.19\textwidth]{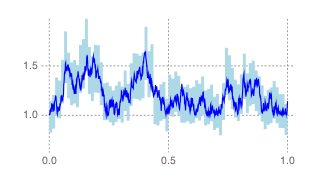} &
\includegraphics[height=0.19\textwidth]{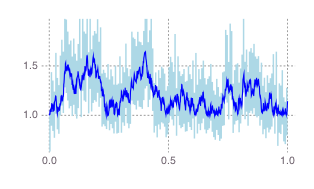} \\
\rotatebox[origin=rt]{-90}{$n=16\,001$\pd}&
\includegraphics[height=0.19\textwidth]{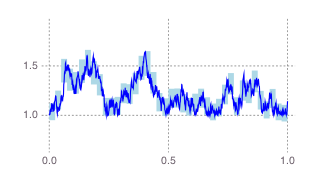}&
\includegraphics[height=0.19\textwidth]{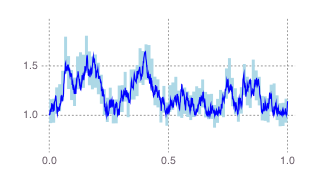}&
\includegraphics[height=0.19\textwidth]{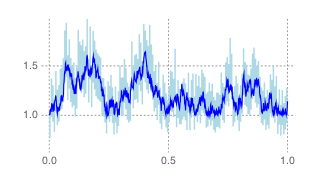}\\
\rotatebox[origin=rt]{-90}{$n=32\,001$\pd}&
\includegraphics[height=0.19\textwidth]{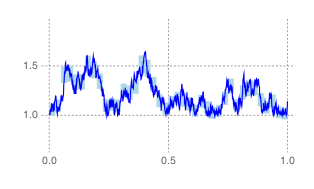}&
\includegraphics[height=0.19\textwidth]{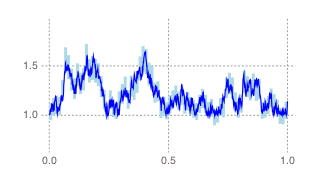}&
\includegraphics[height=0.19\textwidth]{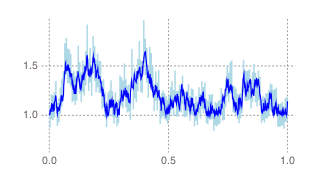}
\end{tabular}
}
\caption{Estimation results for $s_2$ with varying number of observations and bins, no drift.
The solid curve is the true dispersion coefficient while the light blue areas are 98\,\% marginal posterior bands. $n$ is the number of observations, $N$ is the number of bins.
}
\label{fig:s4}
\end{figure}		

\subsection{Inverse Gamma Chain Markov prior}\label{section:improved}

In this section we postulate another prior on the values of the bins. Rather than assuming independent inverse Gamma distributions as in Section~\ref{section:lpp}, we take an \emph{inverse Gamma Markov chain (IGMC)} prior. Posterior inference is still straightforward to implement via Gibbs sampling, as the full conditional distributions are available explicitly and turn out to be inverse Gamma again. We will compare the performance of under both priors.
The results of this section are based on \cite{gugu18b}.

We still consider the diffusion process that solves the SDE \eqref{sde}, with the observations $\scr{X}_n =\{X_{t_{0,n}}, \ldots, X_{t_{n,n}}\}$. We also keep  on taking the drift equal to zero in our statistical analysis. A main reason to depart from the independent prior setting is that this has only limited possibility for adaptation to the local structure of the volatility coefficient, which may become an issue if the volatility has a wildly varying curvature on the time interval $[0,1]$. 

The alternative, fundamentally different, approach, we undertake here, is inspired by ideas in \cite{cemgil07} in the context of audio signal modelling, see also \cite{dikmen10} and \cite{peeling08}, different from the SDE setting that we consider. Instead of using a prior on the (squared) volatility that has piecewise constant realisations on $[0,1]$ with independent coefficients $\theta_k$'s, as in Section~\ref{section:lpp}, we will assume that the sequence $\{\theta_k\}$ forms a suitably defined Markov chain. An immediately apparent advantage of using such an approach is that it induces extra smoothing via dependence in prior realisations of the volatility function across different bins.  Specifically, we proceed as follows, first specifying the piecewise constant prior on the levels on the bins, now being a Markov chain. The bins $B_k$ are the same as before.
We again model $s$ as piecewise constant on bins $B_k$, thus
$
s=\sum_{k=1}^{N} \xi_k \ind_{B_k},
$
and with $\theta_k=\xi_k^2$ it follows that 
$s^2=\sum_{k=1}^{N} \theta_k \ind_{B_k}$.

Introduce auxiliary variables $\zeta_k,$ $k=2,\ldots,N$, and define a Markov chain using the time ordering $\theta_1,\zeta_2,\theta_2,  \ldots,\zeta_k,\theta_k,\ldots,\zeta_N,\theta_N$. Transition distributions of this chain are defined as follows: fix hyperparameters $\alpha_1$, $\alpha_{\zeta}$ and $\alpha$, and set
\begin{equation}
\label{formula:prior0}
\theta_1 \sim \ig(\alpha_{1},\alpha_{1}), \quad \zeta_{k+1} | \theta_k \sim \ig(\alpha_{\zeta},\alpha_{\zeta} \theta_k^{-1}), \quad \theta_{k+1} | \zeta_{k+1} \sim \ig(\alpha,\alpha \zeta_{k+1}^{-1}).
\end{equation}
The prior is referred to as an inverse Gamma Markov chain (IGMC) prior, which reflects the fact that the distributions in \eqref{formula:prior0} are inverse Gamma. In fact, it can even be verified by direct computations employing \eqref{formula:prior0} that the full conditional distributions of $\theta_k$'s and $\zeta_k$'s are inverse gamma too, see \cite{gugu18b} for details, which facilitates Gibbs sampling.
As observed in \cite{cemgil07}, there are various ways of defining an inverse Gamma Markov chain. The point to be kept in mind is that the resulting posterior should be computationally tractable, and the prior on $\theta_k$'s should have a capability of producing realisations with positive correlation structures, as this introduces smoothing among the $\theta_k$'s in adjacent bins. The definition of the IGMC prior in the present work, that employs the latent variables $\zeta_k$'s, takes care of both these important points. For an additional discussion see \cite{cemgil07}. 

In the practical examples presented below, both the number of bins, $N$, and the parametes $\alpha, \alpha_1, \alpha_\zeta$ are considered as hyperparameters. Assigning values to them can be done in multiple ways, see \cite{gugu18b}, but it has been observed in our simulation examples that inferential conclusions with the IGMC prior are quite robust with respect to the choice of $N$. This is because  through the hyperparameters $\alpha$ and $\alpha_{\zeta}$, the IGMC prior has an additional layer for controlling the amount of applied smoothness; when $\alpha$ and $\alpha_{\zeta}$ are equipped with a prior (as we do), they can in fact be learned from the data.

\begin{figure}
\begin{center}
\includegraphics[width=0.48\textwidth]{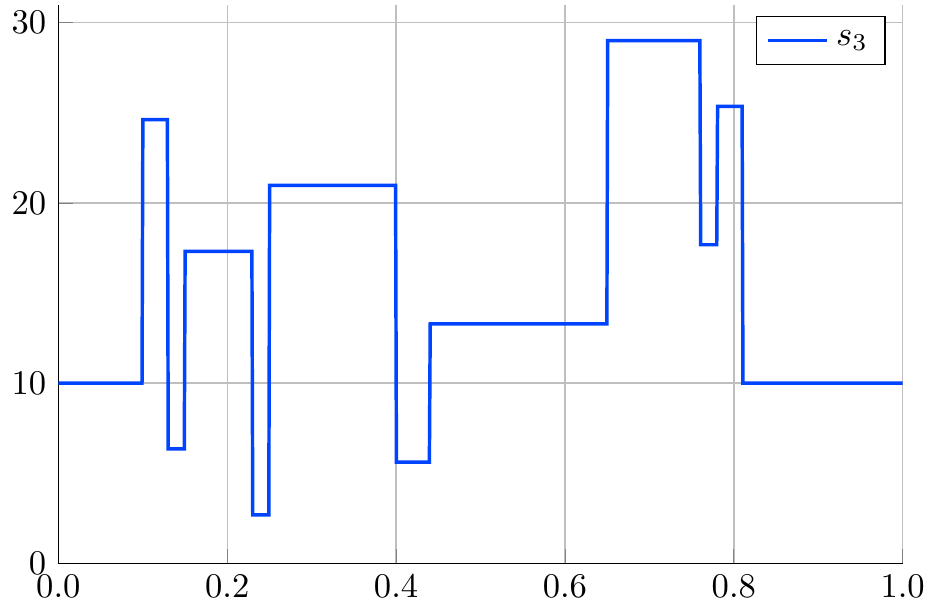}
\caption{The blocks function $s$ from Equation \eqref{s3}.}
\label{fig:blockspath}
\end{center}
\end{figure}

The test function to be estimated is the blocks function from \cite{donoho95} visualized in Figure \ref{fig:blockspath}. With a vertical shift for positivity, this is defined as follows:
\begin{equation}\label{s3}
s(t)= 10 + 3.655606\times\sum_{j=1}^{11} h_j K(t-t_j), \quad t\in[0,1],
\end{equation}
where $K(t)=(1+\operatorname{sgn}(t))/2$, and
\begin{align*}
\{t_j\}&=(0.1,0.13,0.15,0.23,0.25,0.4,0.44,0.65,0.76,0.78,0.81),\\
\{h_j\}&=(4,-5,3,-4,5,-4.2,2.1,4.3,-3.1,2.1,-4.2).
\end{align*}
Our main goal here was to compare the performance of the IGMC prior-based approach to the IIG prior-based one from Section~\ref{section:lpp}. To complete the SDE specification, the true drift coefficient was chosen to be a rather strong linear drift $a_2(x) = -10x+20.$

\begin{figure}
\includegraphics[width=0.48\textwidth]{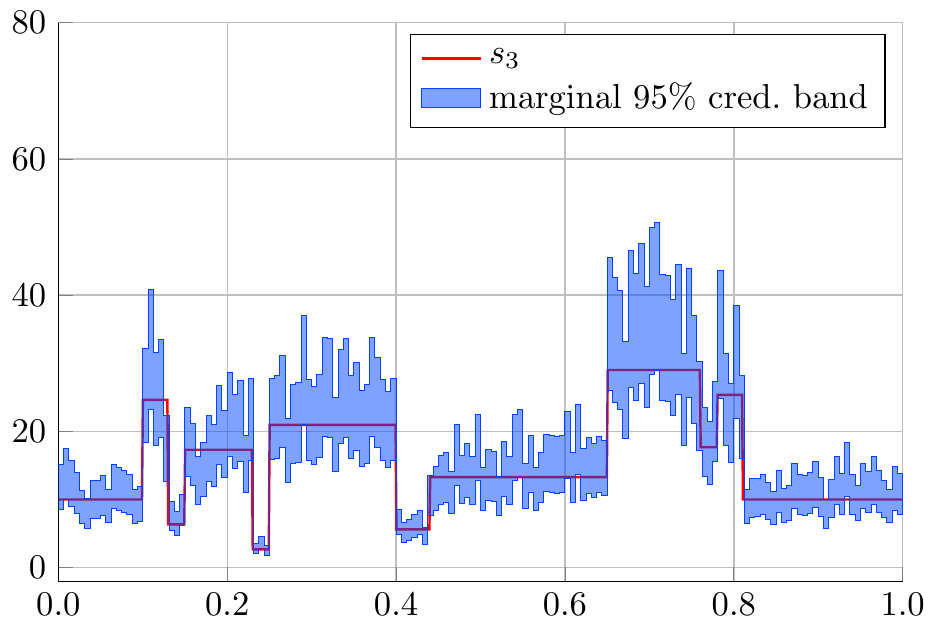}
\includegraphics[width=0.48\textwidth]{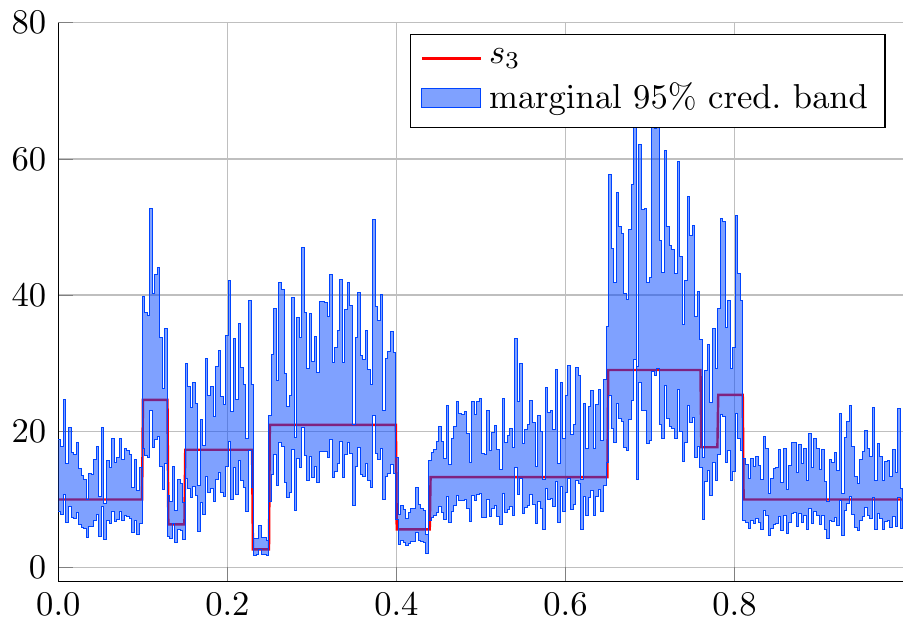}

\bigskip

\includegraphics[width=0.48\textwidth]{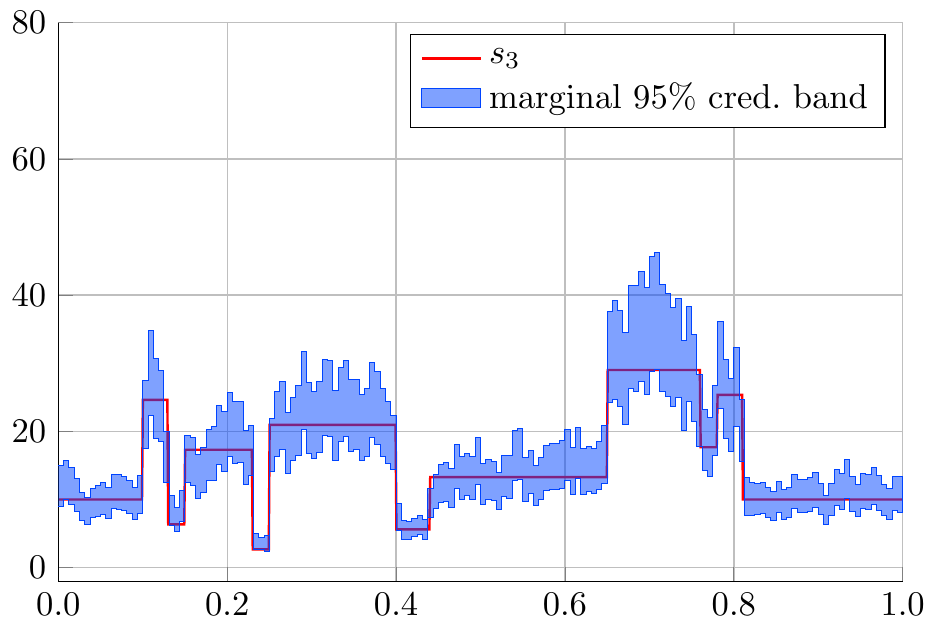}
\includegraphics[width=0.48\textwidth]{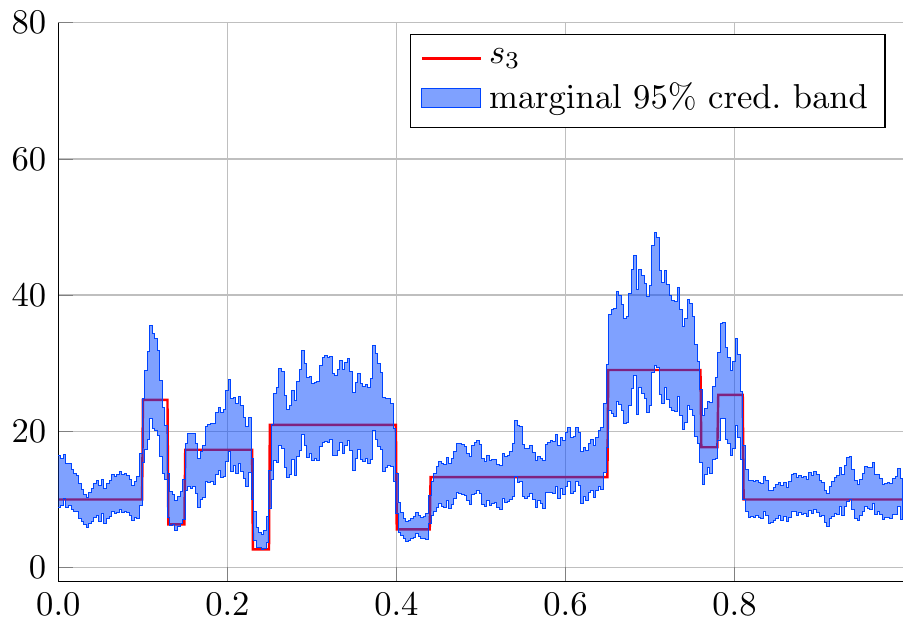}
\caption{Volatility function $s$ with superimposed $95\%$ marginal credible band for the IIG prior $\ig(0.1,0.1)$, using $N=160$ (top left panel) and $N=320$ bins (top right panel). Volatility function $s$ from \eqref{s3} with superimposed $95\%$ marginal credible band for the IGMC prior, using $N=160$ (bottom left panel) and  $N=320$ bins (bottom right panel); in both cases, $\alpha_1=0.1$ and $\alpha =\alpha_{\zeta} \sim\ig(0.3,0.3)$.}
\label{fig:blocks:comp}
\end{figure}

The top left and top right panels of Figure~\ref{fig:blocks:comp} give estimation results obtained with the IIG prior-based approach from Section~\ref{section:lpp}. The number of bins was $N=160$ or $N=320$, and in both these cases we used diffuse independent $\ig(0.1,0.1)$ priors on the coefficients of the (squared) volatility function (see \cite{gugu17} for details). These results have to be contrasted to those obtained with the IGMC prior, plotted in the bottom left and bottom right panels of Figure~\ref{fig:blocks:comp}, where we assumed $\alpha_1=0.1$ and $\alpha =\alpha_{\zeta} \sim\ig(0.3,0.3)$. The following conclusions emerge from Figure~\ref{fig:blocks:comp}:
\begin{itemize}
\item Although both the IGMC and IIG approaches recover globally the shape of the volatility function, the IIG approach results in much greater uncertainty in inferential conclusions, as reflected in wider marginal confidence bands. The effect is especially pronounced in the case $N=320$, where the width of the band for the IIG prior renders it almost useless for inference.
\item The bands based on the IGMC prior look more `regular' than the ones for the IIG prior.
\item The method based on the IIG prior is sensitive to the bin number selection: compare the top left panel of Figure~\ref{fig:blocks:comp} using $N=160$ bins to the top right panel using $N=320$ bins, where the credible band is much wider. On the other hand, the method based on the IGMC prior automatically rebalances the amount of smoothing it uses with different numbers of bins $N$, thanks to the hyperprior on the parameters $\alpha,\alpha_{\zeta}$; in fact, the bottom two plots in Figure~\ref{fig:blocks:comp} look similar to each other.
\end{itemize}

We have shown more pronounced results when using the IGMC prior than by using the IIG prior, although the latter already shows optimal asymptotic properties, as reported in Section~\ref{section:asympvol}. We will see similar behaviour in Section~\ref{sec:poisson_numerical}.

\subsection{Closing remarks}

We have studied nonparametric Bayesian estimation of a volatility function, and presented posterior contraction results when the prior was of independent inverse gamma type. We have given numerical results and have also illustrated that the performance of the statistical procedure had increased when switching from independent priors on the coefficients to a Markov chain prior. The theoretical results can also be applied to real data, see \cite{gugu17} for an example with exchange rate data   and Section~\ref{section:micro} for a follow up with partial observations.

\section{Volatility estimation in presence of microstructure noise}\label{section:micro}

The goal of this section is nonparametric Bayesian estimation of the volatility in presence of market microstructure noise. The exposition is based on \cite{gugushvili2022}, which in turn uses some ideas from \cite{gugu18b}, whose findings have also been used in Section~\ref{section:improved}.

\subsection{Problem formulation}

As in Section~\ref{section:diffusion}, let $X$ solve the one-dimensional stochastic differential equation (SDE) \eqref{sde}. Assume (not necessarily uniformly spaced) observation times $\{t_i: 1\le i\le n\}$ and suppose our observations are
$
\mathcal{Y}_n =\{Y_1, \ldots, Y_{n}\}.
$
Here
\begin{equation}
	\label{eq:obs}
	Y_{i}=X_{t_{i}}+V_{i}, \quad 0<t_1<\cdots<t_n=1,
\end{equation}
and $\{ V_i \}$ is a sequence of independent and identically distributed random variables, independent of the Wiener process $W$ driving \eqref{sde}. The measurement errors $\{V_{i}\}$ are referred to as  microstructure noise. They model such features of financial time series as their discreteness or approximations due to market friction. For widely spaced data in time, microstructure noise can be ignored without great consequences for inferring the volatility function. However, for densely spaced data in time (the so-called high-frequency data regime), this ceases to be true and even small amount of the microstructure noise adversely affects the volatility estimates.

As referenced in \cite{gugushvili2022}, there exists a large body of literature on nonparametric volatility estimation under  microstructure noise. The existing work mainly focusses on estimation of the integrated volatility and is frequentist in its nature. Here we take a Bayesian approach instead. General advantages of the latter have been sufficiently highlighted in the preceding sections.

The key steps of our approach are as follows: we  first reduce our model to the Gaussian linear state space model, which is computationally convenient. Next we specify an inverse Gamma Markov chain prior on the volatility function $s$. This allows us to employ for posterior inference the Gibbs sampler. The latter incorporates a classic Forward Filtering Backward Simulation (FFBS) step.

After detailing the algorithm, we demonstrate the performance of our method on simulated and real data examples, and also provide a theoretical result on the limiting behaviour of the Gamma Markov chain prior.

\subsection{Linear state space model}

We rigorously derive our Bayesian algorithm in the case when the ``true'', data-generating volatility $s$ is a deterministic function of time $t$.  Next, if the ``true'' $s$ is in fact a stochastic process, we simply apply the algorithm without further changes, as if $s$ were deterministic.  The intuition is as follows: first the stochastic volatility is generated, and given a realisation of it, the observations are generated by an independent Brownian motion $W$. 

Denote $t_0=0$. From equation \eqref{sde},
\begin{equation}
	\label{eq:sde}
	X_{t_{i}}=X_{t_{i-1}}+\int_{t_{i-1}}^{t_{i}}a(t,X_t)\dd t+\int_{t_{i-1}}^{t_{i}} s(t)\dd W_t.
\end{equation}

Due to the properties of the Wiener process paths, over short time intervals $[{t_{i-1}},{t_{i}}]$ the term $\int_{t_{i-1}}^{t_{i}} s(t)\dd W_t$ in \eqref{eq:sde} dominates the term  $\int_{t_{i-1}}^{t_{i}}a(t,X_t)\dd t$. Indeed, the latter scales as $\Delta t_i$, whereas the former as $\sqrt{\Delta t_i}$. Following \cite{gugu17}, as in Section~\ref{section:diffusion}, we act as if the process $X$ had a zero drift, $a\equiv 0$.

We thus \emph{assume}
$
X_{t_{i}}=X_{t_{i-1}}+U_{i},
$
where 
$
U_{i}=\int_{t_{i-1}}^{t_{i}} s(t)\dd W_t.
$
This implies that
\begin{equation}
	\label{ui}
	U_{i} \sim N(0,w_i)\quad \text{with} \quad w_i=\int_{t_{i-1}}^{t_{i}} s^2(t)\dd t.
\end{equation}
Furthermore, $\{U_{i}\}$ is a  sequence of independent random variables, for so are the increments of a Wiener process over disjoint time intervals. To lighten the notation, in the sequel we will write $x_i=X_{t_{i}}$, $y_i=Y_{i}$, $u_i=U_{i}$, $v_i=V_{i}$. Our model reduces to the \emph{linear state space model}
\begin{equation}
	\label{eq:ssm}
	\begin{split}
		x_i&=x_{i-1}+u_i,\\
		y_i&=x_i+v_i,
	\end{split}
\end{equation}
where $i = 1, \dots, n$.
In the time series literature, the first equation in \eqref{eq:ssm} is called the \emph{state equation}, whereas the second equation is referred to as the \emph{observation equation}. We assume that $\{v_i\}$ is a sequence of independent $N(0,\eta_v)$ distributed random variables, also independent of the Wiener process $W$ in \eqref{sde}. As a result, $\{v_i\}$ is independent of $\{u_i\}$. We equip the \emph{initial state} $x_0$ with the $N(\mu_0,C_0)$ prior distribution. Then \eqref{eq:ssm} turns into a \emph{Gaussian} linear state space model. The latter is most convenient computationally, as one does not have to deal with an \emph{intractable likelihood}. The latter constitutes the main \emph{computational bottleneck} for Bayesian inference in SDE models.

\subsection{Prior}\label{section:microig}

For the measurement error variance $\eta_v$, we assume a priori an inverse Gamma distribution, $\eta_v \sim \ig(\alpha_v,\beta_v)$. This choice is dictated by conjugacy considerations. In general, the inverse gamma prior on the variance parameter in Gaussian models is  common practice in Bayesian statistics.

The construction of the prior for $s$ is more involved. Fix an integer $m<n$. Then we can write $n=mN+r$  with $0\leq r<m$, where $N=\lfloor {n}/{m}\rfloor$. Define bins $B_k=[t_{m(k-1)},t_{mk})$, $k=1,\ldots,N-1$, and $B_N=[t_{m(N-1)},1]$. Assume a priori
\begin{equation}\label{eq:prior-s}
	s=\sum_{k=1}^{N} \xi_k \ind_{B_k}, 
\end{equation}
where  $N$ (the number of bins) is a hyperparameter. Then $
s^2=\sum_{k=1}^{N} \theta_k \ind_{B_k},
$
where  $\theta_k=\xi_k^2$. 
The prior specification for  $s$ is completed by  assigning a prior distribution to the coefficients $\theta_{1:N}$ as in Section~\ref{section:improved}, which we repeat here. We introduce auxiliary variables $\zeta_{2:N}$, and suppose the sequence $\theta_1,\zeta_2,\theta_2,\ldots,\zeta_k,\theta_k,\ldots,\zeta_N,\theta_N$ forms a Markov chain (in this order of variables). The transition distributions of the chain are defined  by
\begin{equation}
\label{formula:prior}
\theta_1 \sim \ig(\alpha_1,\beta_1), \quad \zeta_{k+1} | \theta_k \sim \ig(\alpha,\alpha \theta_k^{-1}),  \quad \theta_{k+1} | \zeta_{k+1} \sim \ig(\alpha,\alpha \zeta_{k+1}^{-1}),
\end{equation}
where $\alpha_1,\beta_1,\alpha$ are hyperparameters. We refer to this chain as an inverse Gamma Markov chain, see \cite{cemgil07}. 
The variables $\theta_1,\ldots,\theta_N$ are positively correlated, which imposes \emph{smoothing} across different bins $B_k$: under the prior, the small or large values of $\theta$'s on adjacent bins tend to go together. At the same time, we attain partial conjugacy in the Gibbs sampler that we derive below. This leads to simple and tractable MCMC inference.

An uninformative choice $\alpha_1,\beta_1\rightarrow 0$ performs well in practice. The parameter $\alpha$ plays the role of the smoothing parameter.  We equip it with a log-normal prior  $\log \alpha\sim N(a,b)$. The hyperparameters $a\in\mathbb{R},b>0$ are  chosen so as to render the hyperprior on $\alpha$ diffuse. The interpretation of hyperparameter $N$ (or equivalently $m$) is somewhat different than merely a smoothing parameter. It can also be considered as the resolution at which one wants to infer the volatility function. Under finite amounts of data, this resolution cannot be made arbitrarily fine. On the other hand, inference with the IGMC prior is robust with respect to a wide range of values of $N$ due to the fact that the corresponding Bayesian procedure has an additional \emph{regularisation parameter} $\alpha$ that is estimated from the data. See also the earlier Section~\ref{section:improved} and Section~\ref{section:poisson} for related phenomena.

\subsection{Gibbs sampler for sampling from the posterior}

An expression for the posterior of $s$ can be written down in closed form, but does not scale up computationally. The problem is circumvented by following a data augmentation approach: we treat  $x_{0:n}$ as  \emph{missing data}, whereas the $y_{1:n}$ are the observed data; cf.~\cite{tanner87}. An expression for the joint density of all the quantities involved in the model is easily obtained as
\begin{align*}
	\lefteqn{p(y_{1:n}, x_{0:n}, \th_{1:N}, \zeta_{2:N}, \alpha, \eta_v)  =} \\
	&\qquad \left(\prod_{k=1}^n p(y_k\mid x_k, \eta_v)\right)  p(x_{0:n} \mid \th_{1:N}) \\
	&\qquad \times p(\th_1)  \prod_{k=1}^{N-1} \left[ p(\zeta_{k+1}\mid \th_k, \alpha) p(\th_{k+1} \mid \zeta_{k+1}, \alpha)\right] p(\alpha) p(\eta_v).
\end{align*}
Except for $p(x_{0:n} \mid \th_{1:N})$, all the densities have been given in the previous subsection. To obtain an expression for the latter, define (with $\Delta_i\equiv \Delta t_i$)
\begin{align*}
	Z_k &= \sum_{i=(k-1)m+1}^{km} \frac{(x_i-x_{i-1})^2}{\Delta_i}, \quad k=1,\ldots,N-1,\\
	Z_N &= \sum_{i=(N-1)m+1}^{n} \frac{(x_i-x_{i-1})^2}{\Delta_i}.
\end{align*}
Next set $m_k=m$ for  $k=1,\ldots,N-1$, and let $m_N=m+r$. It follows that 
\[ p(x_{0:N} \mid  \th_{1:N}) \propto e^{(x_0-\mu_0)^2/(2C_0)} \prod_{k=1}^N  \th_k^{-m_k/2} \exp\left(-\frac{ Z_k}{2\th_k} \right). \]

Now we can use the Gibbs sampler to sample from the joint conditional distribution of $(x_{0:n}, \th_{1:N}, \zeta_{2:N}, \eta_v, \alpha)$ given observations $y_{1:n}$. The full conditionals of $\th_{1:N}$, $\zeta_{2:N}$, $\eta_v$  required in the Gibbs sampler are easy to derive and are of the inverse Gamma type. Furthermore, the parameter $\alpha$ is updated via a Metropolis-Hastings step; see \cite{gugushvili2022} for all the details. This is quite similar conceptually to Section~\ref{sec:gamma} below. On the other hand, for updating $x_{0:N}$, conditional on all other quantities in the model, we  use the classical Forward Filtering Backward Simulation (FFBS) algorithm for Gaussian state space models (cf.\ Section 4.4.3 in \cite{petris09}). The Gibbs sampler that results from these steps is summarized in Algorithm~\ref{pseudocode}. For details, see \cite{gugushvili2022}.

\begin{algorithm}
	\SetAlgoLined
	\KwData{Observations $y_{1:n}$}
	Hyperparameters $\alpha_1,$ $\beta_1$, $\alpha_v$, $\beta_v$, $a$, $b$, $N$\;
	\KwResult{Posterior samples $\theta_{1:N}^{i}:i=1,\ldots,M$}
	Initialization $\theta_{1:N}^0,$ $\zeta_{1:N}^0$, $\eta_v^0$, $\alpha^0$\;
	\While{$i \leq M$}{
		sample $x_{0:n}^i$ via FFBS\;
		sample $\theta_{1:N}^i$ from the inverse Gamma full conditionals\;
		sample $\zeta_{2:N}^i$ from the inverse Gamma full conditionals\;
		sample $\eta_v^i$ from the inverse Gamma full conditional\;
		sample $\alpha^i$ via a Metropolis-Hastings step\;
		set $i=i+1$. 
	}
	\caption{Gibbs sampler for volatility learning}
	\label{pseudocode}
\end{algorithm}

\subsection{Asymptotics}

In this section we provide some results on the limit behaviour of the Gamma Markov Chain prior. They shed additional light on the prior and are expected to be of use when establishing frequentist asymptotic properties of our nonparametric Bayesian volatility estimation procedure. 

We start with results on the prior conditional mean and variance of the $\theta_k$. These results will be exploited to find an asymptotic regime for the $\theta_k$ when the number of bins tends to infinity. To simplify the exposition, in this section we depart from the setting with $\alpha$ as a random hyperparameter, but instead we assume it to be a deterministic sequence. 
The proof of the following proposition is in \cite{gugushvili2022}.

\begin{proposition}\label{prop:tmv}
	If $\alpha>2$ is a fixed parameter then 
	the IGMC prior satisfies
	\begin{align}
		\EE_k [\theta_{k+1}-\theta_k] 
		& =  \frac{1}{\alpha-1}\theta_k, \label{eq:cmean} \\
		\VV_k (\theta_{k+1}-\theta_k) 
		& =  \frac{\alpha(2\alpha-1)}{(\alpha-1)^2(\alpha-2)} \theta_k^2, \label{eq:cvar}
	\end{align}	
	where $\EE_k$ and $\VV_k$ respectively denote expectation and  variance, conditional on $\theta_k$. Hence, the conditional mean squared error $\EE_k(\theta_{k+1}-\theta_k)^2$ equals $\frac{2(\alpha+1)\theta_k^2}{(\alpha-1)(\alpha-2)}$.
\end{proposition}


As seen from Proposition~\ref{prop:tmv}, the conditional mean squared error in Proposition~\ref{prop:tmv} decreases in $\alpha$. This, in turn, illustrates the regularising property of this parameter.

We are interested in the  behaviour of the prior distribution on the $\theta_k$ for large values of $\alpha$. We will suitably scale $\alpha$ with the number of bins $N$ to obtain a limit result by applying Donsker's theorem (cf.\ Corollary VII.3.11 in \cite{jacod2013limit}). 
Take $\alpha= \gamma N$, where $\gamma$ is a positive scaling factor. We are thus interested in the law of $\theta_k$ for $N\to\infty$. This entails simultaneously increasing the number of bins, as well as the dependence of the values on the bins.  Then, as shown in \cite{gugushvili2022}, with $\theta_1$ fixed, the process $
t\mapsto s^2(t)$, with $s^2(t)=\sum_{k=1}^{N} \theta_k \ind_{B_k}(t)$, converges weakly to the continuous time process $t\mapsto \theta_1 Z_t$, where 
\begin{equation}\label{eq:limitZ} 
	\log Z_t=\sqrt{\frac{2}{\gamma}} W_t-\frac{1}{\gamma} t.
\end{equation}
As a result, the continuous time approximation $Z$ of the process $\theta_k/\theta_1$, with $Z$ as in \eqref{eq:limitZ}, is the Dol\'eans exponential $\mathcal{E}(\sqrt{\frac{2}{\gamma}} W)$, which satisfies the stochastic differential equation 
\begin{equation}\label{eq:sdez}
	\dd Z_t=\sqrt{\frac{2}{\gamma}}Z_t\,\dd W_t.
\end{equation}

\subsection{Numerical results for the Heston model}

The Heston model (see~\cite{heston93}) is  a popular stochastic volatility model. According to Heston's model, the asset price process, denoted by $S$, is governed by  the SDE
\[
\dd S_t = \mu S_t \dd t + \sqrt{Z_t } S_t \dd W_t,
\]
where the process $Z$ follows the CIR or square root process (see \cite{cox85}),
\begin{equation}\label{eq:ZCIR}
	\dd Z_t = \kappa ( \theta - Z_t ) \dd t + \sigma \sqrt{Z_t} \dd B_t.
\end{equation}
The driving Wiener processes $W$ and $B$ are correlated with correlation $\rho$. A common approach in quantitative finance is to log transform $S$ to $X_t=\log S_t$. Then $X$ obeys a diffusion equation with volatility $s(t)=\sqrt{Z_t}$,
\[
\dd X_t = (\mu-\frac{1}{2} Z_t)\dd t+\sqrt{Z_t}\dd W_t,
\]
which can verified by It\^{o}'s formula.

Assume high-frequency observations on  the log-price process $X$ with additive noise $V_i \sim N(0,\eta_v)$ are available. Note that the volatility $s$ is a random function. No further knowledge of the data generation mechanism is assumed. The Heston model setting is challenging for any nonparametric volatility estimation method and as such serves as a good test example. To generate synthetic data, we used the parameter values
$\mu = 0.05,$
$\kappa = 7,$
$\theta = 0.04,$
$\sigma = 0.6,$
$\rho = -0.6$.
These are similar to the ones obtained via fitting the Heston model to real data  (see Table 5.2 in \cite{vdp06}). Furthermore, the noise variance was $\eta_v = 10^{-6}.$ This guarantees a large enough signal-to-noise ratio in the model \eqref{eq:ssm}, which is a prerequisite for successful estimation of volatility. Finally, the parameter choice $\mu=0.04$ corresponds to  a realistic $4\%$ log-return rate.

\begin{figure}
	\begin{center}
		\includegraphics[width=0.75\textwidth]{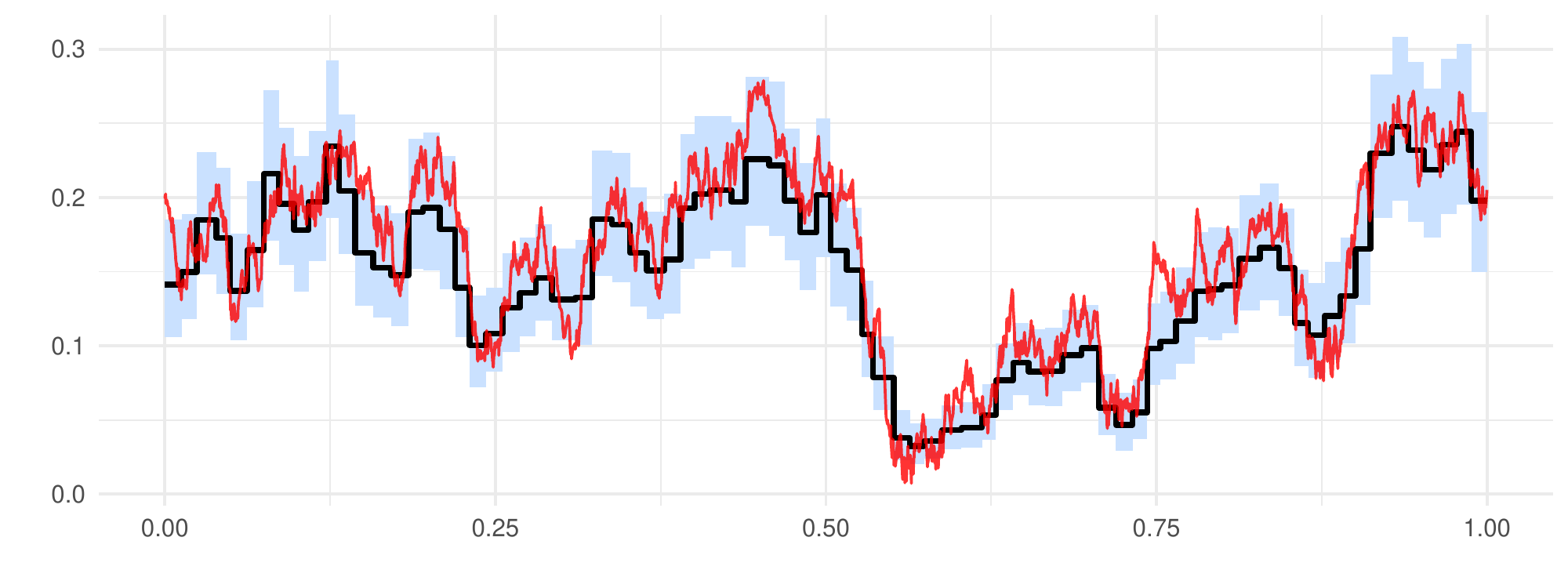}
		\caption{Posterior mean and pointwise $95\%$ credible band for the Heston model.  The true volatility function is plotted in red, the black step function gives the posterior mean, and the credible band is shaded in blue.}
		\label{fig:heston}
	\end{center}
\end{figure}

Inference results with $N=80$ bins are given in Figure~\ref{fig:heston}. These are surprisingly accurate, despite the general difficulty of the problem and the amount of assumptions that went into our estimation procedure. Indeed, the Heston model does not formally satisfy assumptions under which our Bayesian method was derived; this concerns specifically the non-deterministic character of the volatility function $s$.

\begin{figure}
	\begin{center}
		\includegraphics[width=0.75\textwidth]{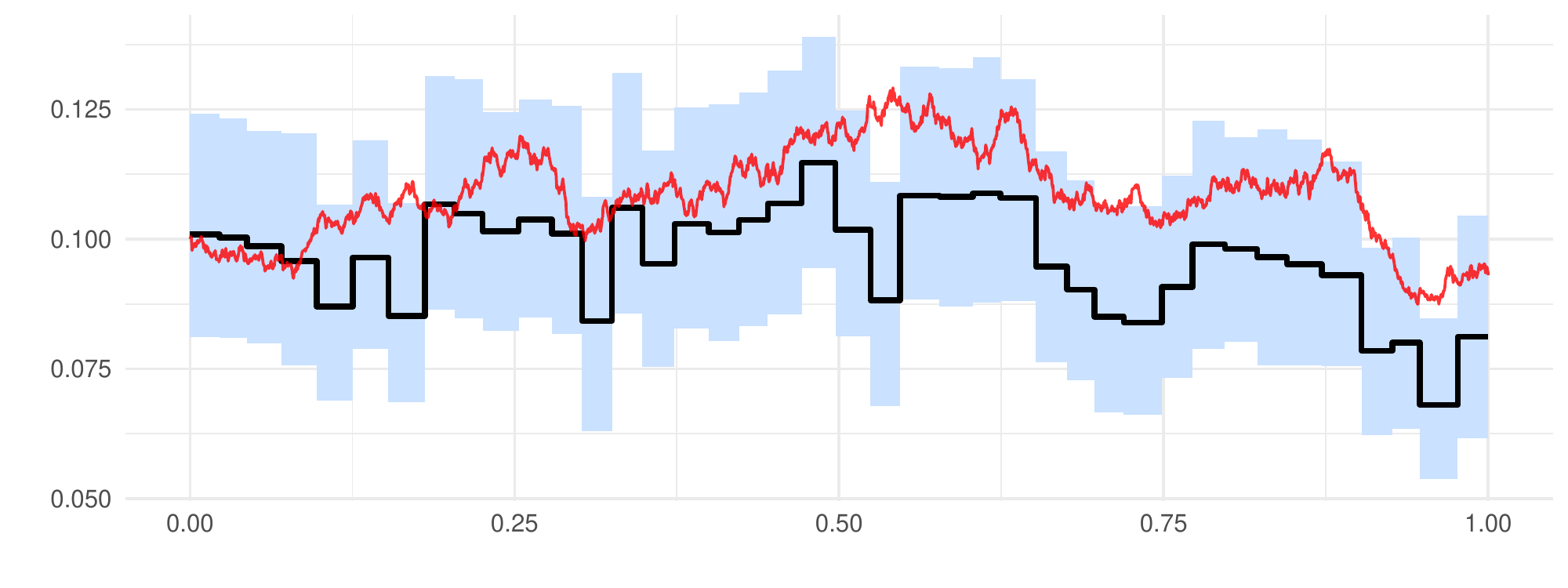}
		\caption{Posterior results for the Heston model where the CIR volatility  is replaced by the root of the continuous time limit of the  prior, $N=40$, starting value of the volatility is $0.1$. The colors have the same meaning as in Figure~\ref{fig:heston}.}
		\label{fig:var_posterior_start40_4}
	\end{center}
\end{figure}

\begin{figure}
	\begin{center}
		\includegraphics[width=0.75\textwidth]{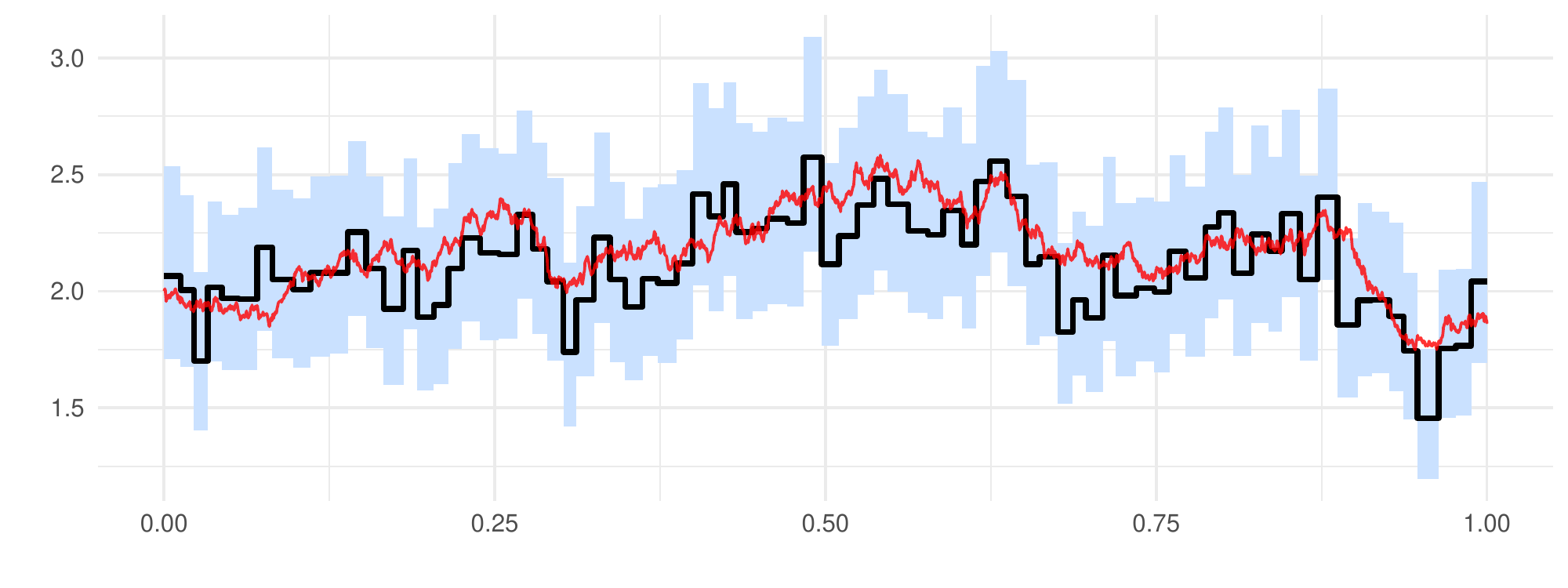}
		\caption{Posterior results for the Heston model where the CIR volatility  is replaced by the root of the continuous time limit of the  prior, $N=80$, starting value of the volatility is $2.0$.}
		\label{fig:var_posterior_start80_4}f
	\end{center}
\end{figure}

We also performed experiments for a variation of the above Heston model, where we replaced the squared volatility process with the limit process $\theta_1 Z$, where  $Z$ is as in \eqref{eq:limitZ}, or, which is the same, as in \eqref{eq:sdez}, with $\sqrt{\frac{2}{\gamma}}=0.6$. These experiments illustrate the influence of the number of bins $N$ and the starting values, see Figures~\ref{fig:var_posterior_start40_4} and \ref{fig:var_posterior_start80_4}. Here the results with a smaller number of bins $N=40$ combined with a low initial value are less satisfactory than with $N=80$ and higher starting value.

\subsection{Exchange rate data}
\label{section:real}
In this section, we infer volatility of the high frequency foreign exchange rate data made available by Pepperstone Limited, the London based forex broker.\footnote{Note that the data are not available from the Pepperstone website any longer, but can be obtained directly from the present authors.} As we shall show shortly, the inferred volatility looks realistic, however there is substantial uncertainty surrounding the inferential results (as demonstrated by the wide marginal credible band).

\begin{figure}
	\begin{center}
		\includegraphics[width=0.75\textwidth]{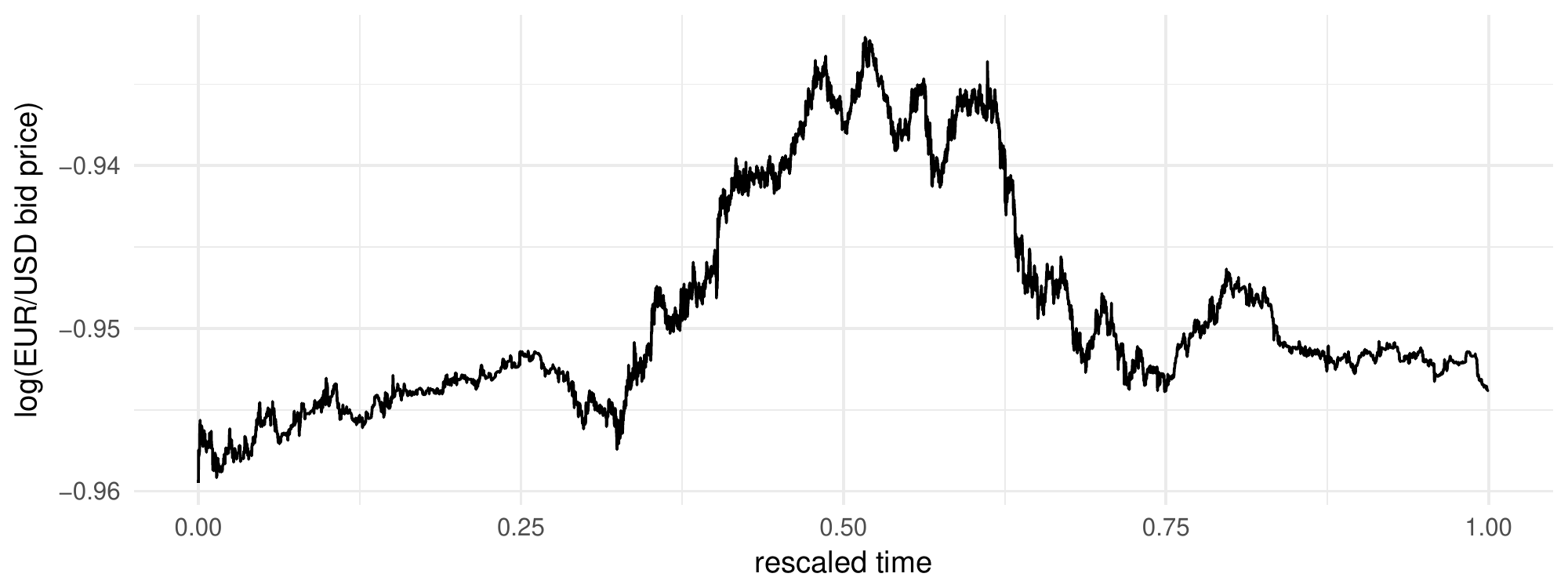}
		\includegraphics[width=0.75\textwidth]{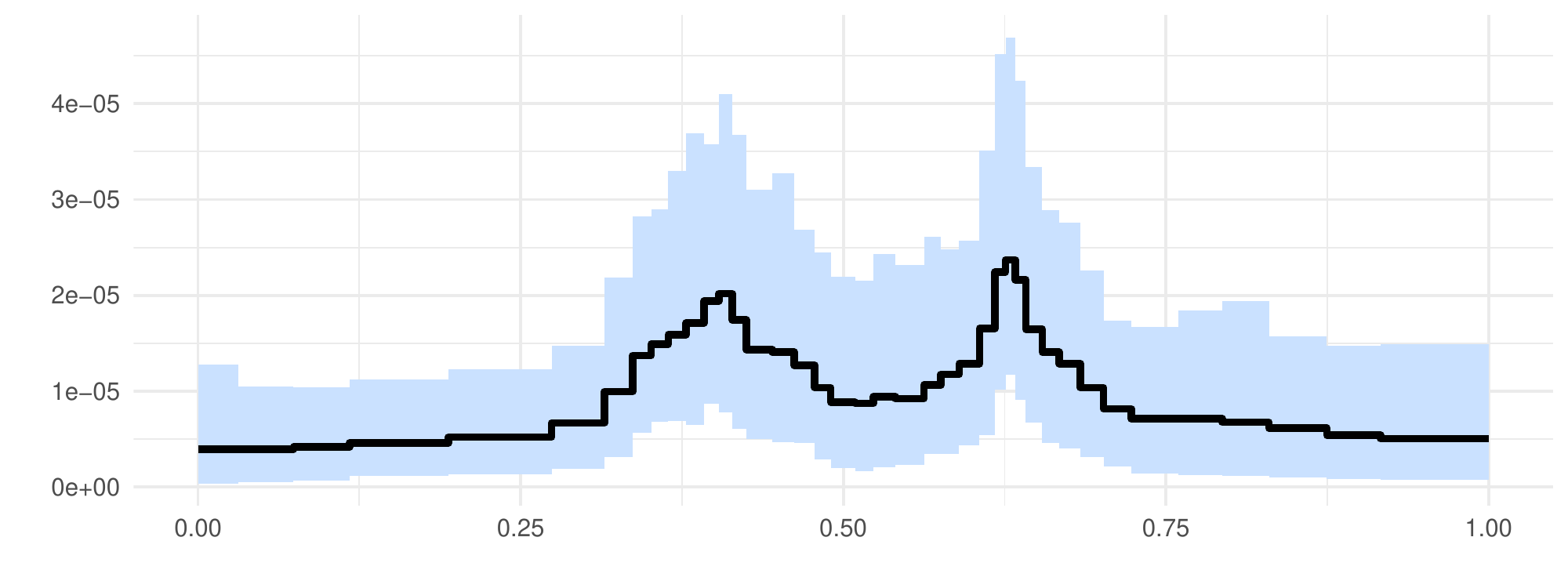}
		\caption{Top: Natural logarithm of the EUR/USD exchange rate data for 2 March 2015 analysed in Section~\ref{section:real}. Bottom: Posterior mean (black curve) and pointwise 95\% credible band (blue band) for the volatility function. The time axis is rescaled to $[0,1]$.}
		\label{fig:eurus:data}
	\end{center}
\end{figure}

In greater detail, we use the EUR/USD tick data (bid prices) for 2 March 2015. Note that the independent additive measurement error model \eqref{eq:obs} becomes harder to justify for highly densely spaced data in time. Therefore we opted to subsample the data, retaining every $10$th observation. Subsampling is a common strategy in this context; see, e.g., Section 5 in \cite{mancini2015spot}. In our case subsampling results in a total of $n=13\, 025$ observations over one day, about~9 per minute. As in the case of Heston model, we infer the volatility from the log of the observed time series. We did not apply any other data preprocessing steps (such as detection and removal of jumps from the observed data), because we wanted to focus on our Bayesian procedure itself.  The data are shown in Figure~\ref{fig:eurus:data}, top panel. Note that time is rescaled to the interval $[0,1]$, so that $t=0$ corresponds to the midnight and $t=0.5$ to  noon.

\begin{figure}
	\begin{center}
		\includegraphics[width=0.75\textwidth]{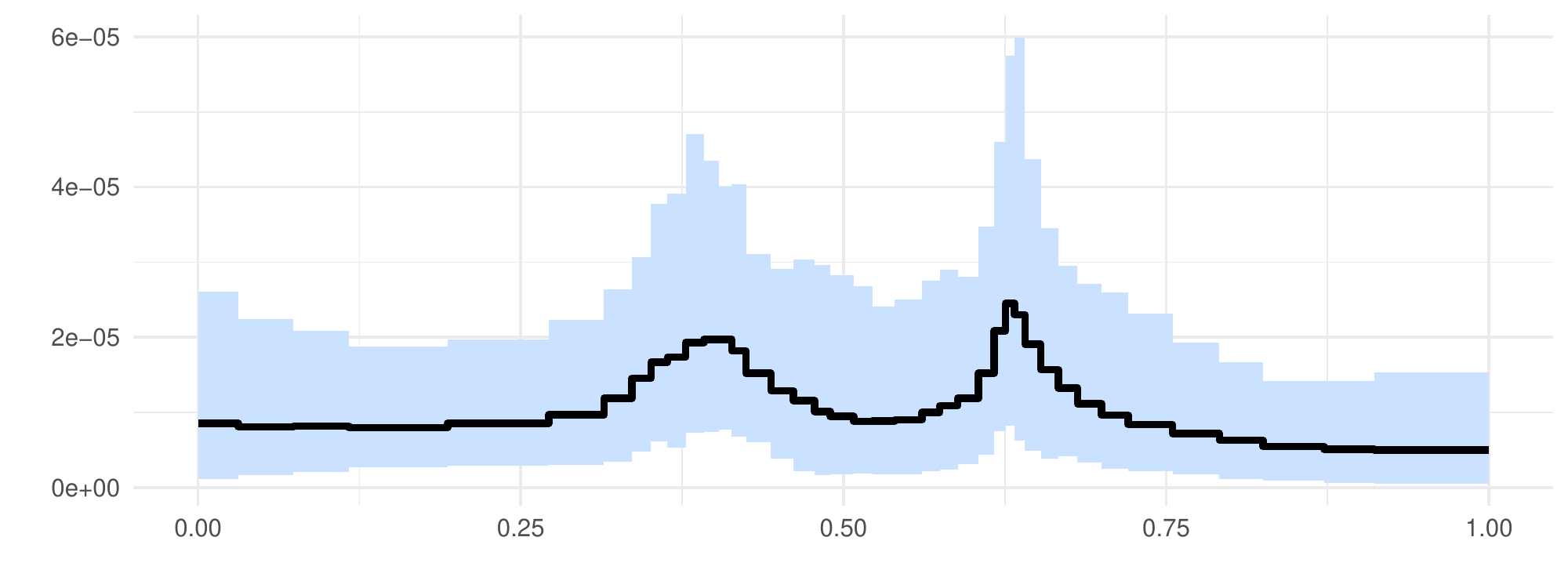}
		\caption{Posterior mean (black curve) and pointwise 95\% credible band (blue band) for the volatility of the further subsampled EUR/USD exchange rate data analysed in Section~\ref{section:real}. The time axis is rescaled to $[0,1]$.}
		\label{fig:eurus2}
	\end{center}
\end{figure}

The inference results are given in the bottom panel of Figure~\ref{fig:eurus:data}.  We used $N=96$ bins, which corresponds to updating volatility every 15 minutes. There is considerable uncertainty in the inferred volatility. Reassuringly, the volatility is lower during the night hours. It also displays two peaks corresponding to the early morning and late afternoon hours of the day. In Figure~\ref{fig:eurus2} we give inference results obtained via further subsampling of the data, retaining $50\%$ of the already subsampled observations. The posterior mean is quite similar to that in Figure~\ref{fig:eurus:data}, whereas the wider credible band reflects greater inferential uncertainty due to a smaller sample size. The stability of the inferred volatility across different observation frequency regimes partially validates the model and the approach that we use.

\subsection{Conclusions}

Volatility estimation under the presence of microstructure noise is much more challenging than the same task under noiseless observations. As such, accurate inference on volatility under  microstructure noise requires large amounts of data. These, however, are readily available in financial applications. In turn, this necessitates finding a volatility estimation method that scales well with data.

The prior that we employ together with and a deliberate, but asymptotically harmless misspecification of the drift (by taking $b\equiv 0$) enable us to combine our earlier work in \cite{gugu18b}, partially reproduced here in Section~\ref{section:improved}, with the FFBS algorithm for Gaussian linear state space models. This results in an intuitive and fast algorithm  to obtain samples from the posterior.

An interesting future research topic in this context is to explicitly account for the possible presence of jumps in financial time series.

\section{Poisson intensity estimation}\label{section:poisson}

This section is based on \cite{gugushvili2020fast}.
Consider a Poisson point process $X$ on an interval $[0,T]$ of the real line $\mathbb{R}$, which we equip with the Borel $\sigma$-field $\mathcal{B}([0,T])$. Such Poisson point processes are also often called non-homogeneous Poisson processes. The process $X$ is a random integer-valued measure on $[0,T]$ (we assume the underlying (complete) probability space $(\Omega,\mathcal{F},\mathbb{Q})$ in the background), such that
\begin{enumerate}
	\item for any disjoint subsets $B_1,B_2,\ldots,B_m\in\mathcal{B}([0,T]),$ the random variables $X(B_1),$ $X(B_2),\ldots,X(B_m)$ are independent, and
	\item for any $B\in\mathcal{B}([0,T]),$ the random variable $X(B)$ is Poisson distributed with parameter $\Lambda(B)$. Here $\Lambda$ is a finite measure on $([0,T],\mathcal{B}([0,T])),$ called the intensity measure of the process $X.$ Moreover, it is assumed that $\Lambda$ admits a density $\lambda$ with respect to the Lebesgue measure on $\mathcal{B}([0,T])$.
\end{enumerate}
Intuitively, the process $X$ can be thought of as random scattering of points in $[0,T],$ where the way the scattering occurs is determined by properties (i)--(ii).

A popular assumption in the statistical literature (see, e.g., \cite{karr86} and \cite{kutoyants98}) is that one observes independent copies $X_1,\ldots,X_n$ of the process $X$.  Based on these observations, we aim to construct an estimator of the intensity function $\lambda$.

\subsection{Likelihood}
By Proposition~6.11 in \cite{karr86} or Theorem 1.3 in \cite{kutoyants98}, the law of $X$ under the parameter value $\lambda$ admits a density $p(\cdot\,;\lambda)$ with respect to the measure induced by a standard Poisson point process of rate $1.$ This density is given by
\begin{equation*}
p(\xi;\lambda)=\exp\left( \int_{[0,T]} \log \lambda(x)\mathrm{d}\xi(x) - \int_{ [0,T]  } (\lambda(x)-1) \mathrm{d}x \right),
\end{equation*}
where $\xi=\sum_{i=1}^m \delta_{x_i}$ is a realisation of $X$ (here $\delta_{x_i}$ denotes the Dirac measure at $x_i$). 
We assume independent observations $X_1,\ldots, X_n$ with the same distribution as $X$ and write 
$ X_j = \sum_{i=1}^{m_j} \delta_{X_{ij}}$ for $j=1,\ldots, n$. 
If we define $X^{(n)}=(X_1,X_2,\ldots,X_n)$, then it follows that the likelihood $L(X^{(n)};\lambda)$  can be written as
\begin{equation}
\label{likelih}
L(X^{(n)};\lambda)=\prod_{j=1}^n \exp\left( \int_{ [0,T] } \log \lambda(x)\mathrm{d}X_{j}(x) - \int_{ [0,T]  } \left(\lambda(x)-1\right) \mathrm{d}x \right).
\end{equation}

\subsection{Prior}

Fix a positive integer $N$. Let $0=b_0<b_1<\cdots<b_{N-1}<b_N=T$ be a grid of points on the interval $\mathcal{X}=[b_0,b_N]$, for instance a uniform grid. Using this grid, define the bins $B_k=[b_{k-1},b_{k}),$ $k=1,\ldots,N-1$, with the last bin $B_N = [b_{N-1},b_N]$. In order to define a prior on $\lambda$, introduce a collection of piecewise constant functions $\lambda$ on bins $B_k$,
\begin{equation}
\label{lambda}
\lambda(x)=\sum_{k=1}^N \psi_k \mathbf{1}_{B_k}(x), \quad x\in[0,T].
\end{equation}
We assume  $\psi_k \iid \operatorname{G}(\alpha,\beta)$ and  define the prior $\Pi_N$ on $\lambda$ to be the law of the random function \eqref{lambda}.
We will refer to the prior $\Pi_N$ as the independent gamma prior. If the grid on $\mathcal{X}=[0,T]$ is taken to be uniform, the prior has three hyperparameters: $\alpha,$ $\beta,$ and $N.$ 
Depending on the amount of available data and the degree of prior knowledge on $\lambda$, the hyperparameters $\alpha,\beta$ can for instance be chosen to render the $\operatorname{G}(\alpha,\beta)$ prior diffuse (non-informative). This corresponds to the case when one has little information on the magnitude and shape of $\lambda$. 
The hyperparameter $N$ can be viewed as a smoothing parameter and tuned using an empirical Bayes procedure.

Using a prior with piecewise constant realisations as in \eqref{lambda} is not unnatural.
In fact, in practical applications it is often the case that Poisson point processes are only partially observed through aggregate counts of points over successive intervals (cf.~\cite{henderson03}). The intensity function cannot be learned beyond the resolution of these intervals, which lends support to using priors with piecewise constant realisations to model the intensity function.

\subsection{Posterior}
Since the function $\lambda$ in our case is parametrised by the coefficients $\psi_1,\ldots,\psi_N$, the posterior distribution of the intensity function $\lambda$ given the data $X^{(n)}$ can be equivalently described through the posterior distribution of $\psi_1,\ldots,\psi_N$ given  $X^{(n)}$. Transition from the prior to the posterior can be thought of as updating our prior opinion on $\lambda$ upon seeing the data $X^{(n)}.$
\begin{lemma}\label{lem:indepposterior}
	\label{lem:gamma}
	Define 
	\begin{equation}\label{eq:def-Hk}	
	H_k = \sum_{j=1}^n \sum_{i=1}^{m_j} \ind_{\{X_{ij} \in B_k\}}
	\end{equation}
	and $\Delta_k= b_k - b_{k-1}$ for $k=1,\ldots,N.$ Then  $\psi_1,\ldots,\psi_N$ are a posteriori independent and 
	\begin{equation}
	\label{post:gamma}
	\psi_k \mid   X^{(n)} \sim \operatorname{G}(\alpha  + H_k, \beta + n\Delta_k ), \quad k=1,\ldots,N.
	\end{equation}	
\end{lemma}
Thus, the posterior for $\lambda$ is known in closed form. The posterior mean is given by
$ \hat{\lambda}(x)=\sum_{k=1}^N \hat{\psi}_k \mathbf{1}_{B_k}(x)$, $x\in [0,T]$,
with
\begin{equation}
\label{psik}
\hat{\psi}_k = \frac{ H_k + \alpha }{ n\Delta_k+\beta }, \quad k=1,\ldots,N.
\end{equation}
Marginal $1-\gamma$ credible bands for $\lambda$ can be obtained by producing $1-\gamma$ credible intervals for $\psi_k,$ using the lower and upper $\gamma/2$-quantiles of the gamma distribution. This gives Bayesian uncertainty quantification for estimation of $\lambda.$

\subsection{Asymptotics}

In this subsection we perform the asymptotic frequentist analysis of the Bayesian procedure we introduced above. This concerns the study of the asymptotic properties of the posterior measure as the sample size $n\rightarrow\infty$.

For simplicity in the proofs, we assume that the grid $\{b_k\}$ defining the bins $B_k$ is uniform: $b_k=Tk/N,$ $k=0,\ldots,N$, so that the bins are of equal width $\Delta_k=\Delta=T/N$.
Our first result shows that the posterior mean $\hat{\lambda}$ is a consistent estimator of $\lambda_0$, and establishes its convergence rate. The expectation $\ee[\,\cdot\,]$ here and in the sequel is always under the law of the observations with the ``true" parameter value $\lambda_0$.

\begin{theorem}
	\label{thm:mean}
Assume	the intensity function $\lambda_0\colon[0,T]\rightarrow (0,\infty)$ is H\"older continuous: there exist constants $L>0$ and $0<h\leq 1$, such that
	\[
	|\lambda_0(x)-\lambda_0(y)| \leq L |x-y|^h, \quad \forall x,y\in [0,T].
	\]
If $N \asymp n^{1/(2h+1)}$,  then
	\[
	\label{mse}
	\ee [\| \hat{\lambda} - \lambda_0  \|_2^2] \lesssim n^{-2h/(2h+1)}.
	\]
\end{theorem}
The right-hand side is the optimal rate for estimating an $h$-H\"older-regular intensity function, see \cite{kutoyants98}.

The next result gives a posterior contraction rate in the $L_2$-metric. Whereas Theorem~\ref{thm:mean} dealt with the `centre' of the posterior distribution, the theorem below deals with the entire posterior distribution.
\begin{theorem}
	\label{thm:post}
	Let the assumptions of Theorem~\ref{thm:mean} hold, and let $\varepsilon_n\asymp n^{-h/(2h+1)}$. Then, for any sequence $M_n\rightarrow \infty$,
	\[
	\ee [\Pi_N(\|\lambda-\lambda_0\|_2 \ge M_n \varepsilon_n \mid X^{(n)})] \rightarrow 0
	\]
	as $n\rightarrow\infty$.
\end{theorem}
The first conclusion that follows from Theorem~\ref{thm:post} is that the proposed Bayesian procedure is consistent: as the sample size $n\rightarrow\infty$, the posterior puts most of its mass on (shrinking) $L_2$-neighbourhoods around the true parameter $\lambda_0$. Furthermore, the rate $\varepsilon_n\asymp n^{-h/(2h+1)}$ in  Theorem~\ref{thm:post} is the optimal posterior contraction rate for $h$-H\"older-regular intensity functions.

\subsection{Gamma Markov chain prior}
\label{sec:gamma}

The assumption that the coefficients $\psi_k$ of the intensity function $\lambda$ are a priori independent and gamma distributed can be relaxed.    
For this we  use ideas that have appeared in various works in the audio signal modelling literature (see, e.g., \cite{cemgil07}, \cite{dikmen10} and \cite{peeling08}). This is similar to the approach taken in Section~\ref{section:improved}. 
Specifically,  we postulate that the $\{\psi_k\}$ form a gamma Markov chain (GMC). This chain is defined as follows: introduce auxiliary variables $\zeta_k,k=2,\ldots,N$, use the time ordering $\psi_1,\zeta_2,\psi_2,\ldots,\zeta_N,\psi_N$, and set
\begin{equation}
\label{formula:priorpsi}
\psi_1 \sim \g(\alpha_{1},\beta_{1}), \quad \zeta_{k} \mid \psi_{k-1} \sim \ig(\alpha_{\zeta},\alpha_{\zeta} \psi_{k-1}), \quad \psi_{k} \mid \zeta_{k} \sim \g\left(\alpha_{\psi},\frac{\alpha_{\psi}}{\zeta_{k}}\right), 
\end{equation}
where $k\in \{2,\ldots, N\}$. 
The name of the chain reflects the fact that its transition distributions are (inverse) gamma. The parameters $\alpha_1,$ $\beta_1$, $\alpha_{\zeta}$ and $\alpha_{\psi}$ are the hyperparameters of the GMC prior. The hyperparameters $\alpha_1,\beta_1$ allow one to `release' the chain at the origin. This is important to avoid possible edge effects in nonparametric estimation due to a strong specification of the prior at the time origin $t=0$. A principal aim in using latent variables $\zeta_k$'s in \eqref{formula:priorpsi} is to attain positive correlation between $\psi_k$'s. In the intensity function modelling context this induces smoothing across different bins.

The posterior distribution with the GMC prior is not available in closed form. However, ``full'' conditional distributions are known in closed form which facilitates Gibbs sampling.
The full conditional distributions of $\psi_k$, $\zeta_k$ are
\begin{align}
\zeta_k \mid \psi_k,\psi_{k-1} &\sim \ig\left(\alpha_{\zeta}+\alpha_{\psi}, {\alpha_{\zeta}}{\psi_{k-1}}+{\alpha_{\psi}}{ \psi_k}\right), \quad k=2,\ldots,N, \label{fullcond2} \\
\psi_k \mid \zeta_{k+1},\zeta_{k},X^{(n)} &\sim \g\left(\alpha_{\psi}+\alpha_{\zeta}+H_k,\frac{\alpha_{\psi}}{\zeta_k}+\frac{\alpha_{\zeta}}{ \zeta_{k+1}} + n\Delta_k\right), \quad k=2,\ldots,N-1, \label{fullcond1} \\
\psi_1 \mid \zeta_2,X^{(n)} & \sim \g\left(\alpha_1+\alpha_{\zeta}+H_1,  \beta_1 + \frac{\alpha_{\zeta}}{ \zeta_{2}} + n\Delta_1 \right), \label{fullcondstart} \\
\psi_N \mid \zeta_N,X^{(n)} & \sim \g\left(\alpha_{\psi}+H_N,\frac{\alpha_{\psi}}{\zeta_N}+n\Delta_N\right). \label{fullcondend}
\end{align}
The Gibbs sampler cycles through formulae \eqref{fullcond2}--\eqref{fullcondend} to generate samples from the posterior distribution of $\{\psi_k\}$, $\{\zeta_k\}$ given the data $X^{(n)}$.  One can initialize the sampler e.g.\ by providing values for $\psi_1,\ldots,\psi_N$.

\subsection{Numerical results}\label{sec:poisson_numerical}

\begin{figure}
	\begin{center}
		\includegraphics[width=0.9\textwidth]{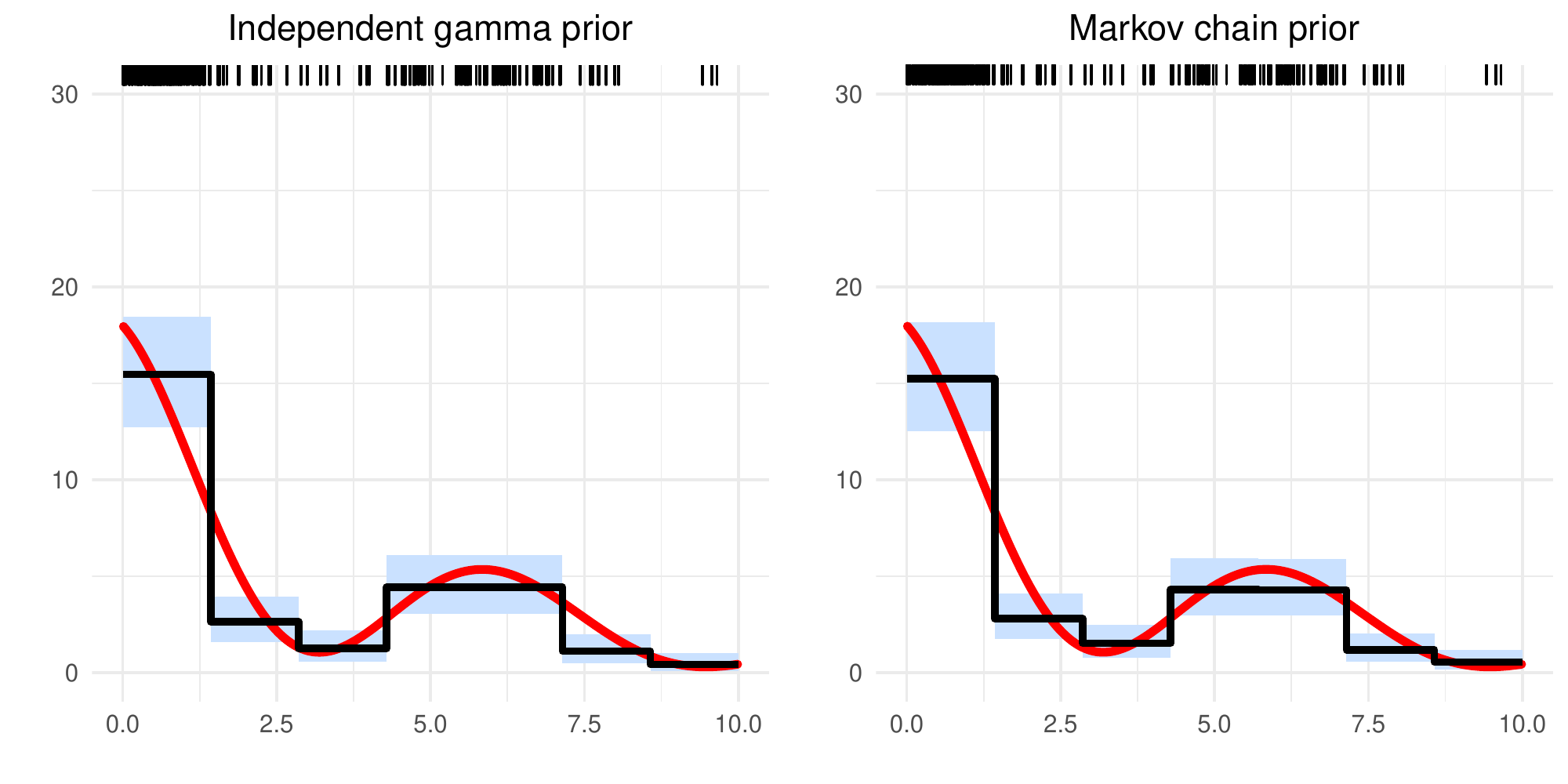}
		\includegraphics[width=0.9\textwidth]{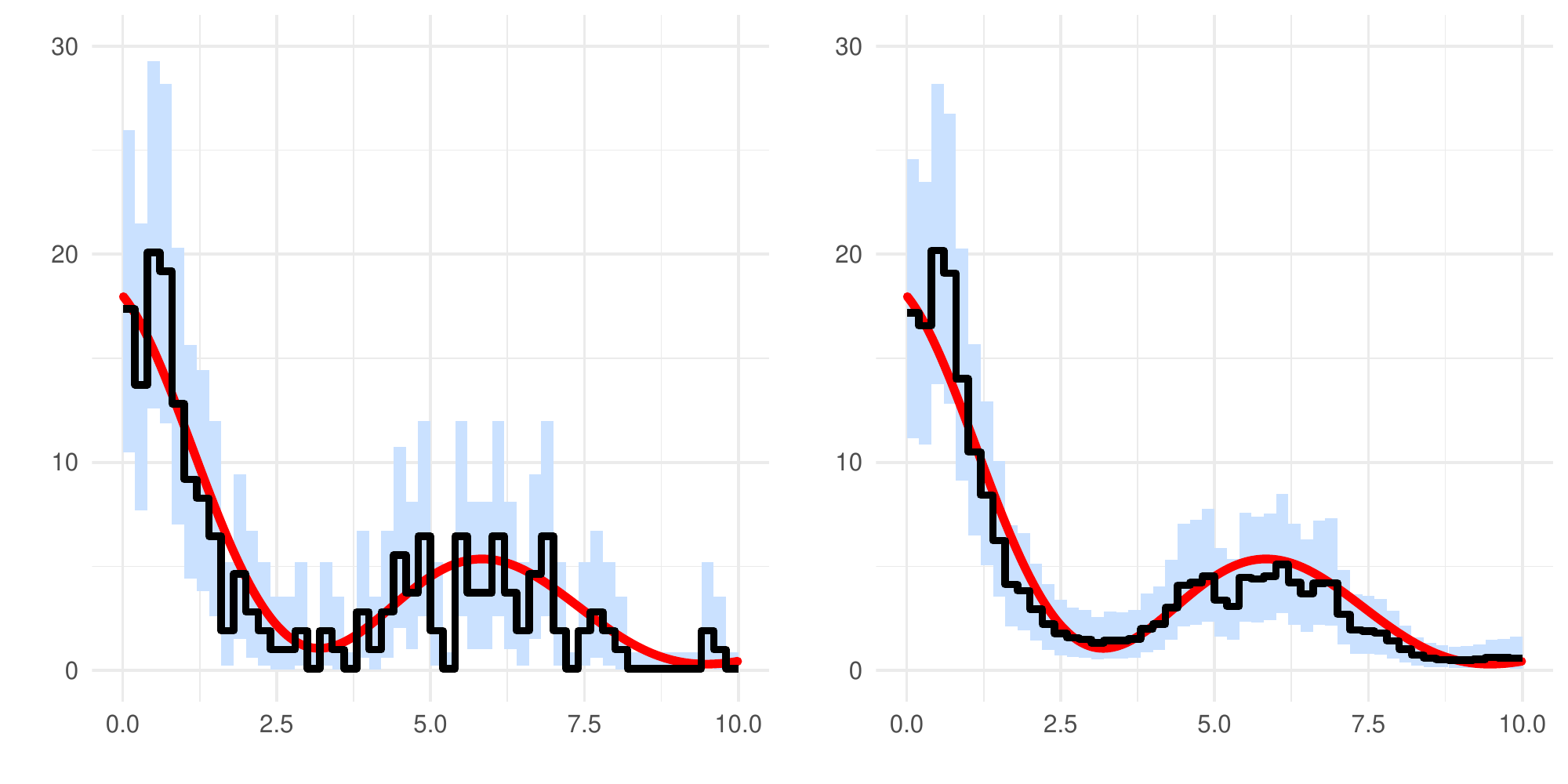}
	\end{center}
	\caption{Estimation results for the oscillating exponential function $\lambda_0(x)=2 e^{-x/5} (5+ 4\cos(x))$ with $n=5$ replicates. The top plot corresponds to the method based on the independent gamma prior, with $N=7$ chosen via the empirical Bayes method (details are in Section 3.3 of \cite{gugushvili2020fast}). The bottom plot corresponds to the method based on the GMC prior with $N=50$ bins.}
	\label{fig:example1n5}
\end{figure}

Consider 
\[
\lambda_0(x)=2 e^{-x/5} (5+ 4\cos(x)) , \quad x\in[0,10].
\]
A principal challenge in inferring this function consists in the fact that it takes small values in the middle part of its domain and has a slope of changing sign. In Figure~\ref{fig:example1n5} the posterior mean and marginal $95\%$ posterior credible bands are shown. From the lower panels one can clearly see the beneficial smoothing effect of the GMC-prior. 

We refer to  \cite{rosas2023bayesian} for an application of the proposed methods in estimating heart rate variability.

\subsection{Closing remarks}
The presented approach is computationally appealing, as it scales with the number of bins and not the sample size. One can think of the chosen number of bins as the resolution at which one estimates the intensity function. If this resolution is too fine, the independent Gamma prior results in  overfitting. This can be avoided by using the GMC-prior, which induces a smoothing effect on neighbouring bins.

\section{Estimation for gamma-driven stochastic differential equations}\label{section:gammasde}

The goal of the present section, based on \cite{bgsslevy}, is Bayesian nonparametric estimation of the positive local scale function $\sigma$ that appears in the L\'evy-driven stochastic differential equation 
\begin{equation} \label{eq:sdegamma}
\dd X_t=\sigma(X_{t-})\,\dd L_t,\,\, X_0= 0, 
\end{equation}
from observations of the solution $X$. Here
$L$ is a gamma process with $L_0=0$ and therefore $L$ is a subordinator, i.e.\ a stochastic process with monotonous sample paths. Furthermore, $L$ has
a L\'evy measure $\nu$ admitting the L\'evy density
\begin{equation}\label{eq:v}
v(x)= \alpha x^{-1}\exp(-\beta x),\,\, x>0,
\end{equation}
where $\alpha$ and $\beta$ are two positive constants. The process $L$ has independent increments,  and $L_t-L_s$ has the Gamma$(\alpha (t-s),\beta)$ distribution for $t>s$, 
i.e.\ the distribution with density $p(x)=\frac{\beta^{\alpha(t-s)}}{\Gamma(\alpha(t-s))}x^{\alpha(t-s)-1}\exp(-\beta x)$, $x\geq 0$.

\subsection{Preliminaries}

Under the assumption that the function $\sigma$ (in view of financial applications we refer to it as volatility function) is measurable and satisfies a linear growth condition, it has been shown in \cite{bgssweak} that Equation~\eqref{eq:sdegamma} admits a weak solution that is unique in law. Note that $X$ is a Markov process. Under the stronger condition that $\sigma$ is 
Lipschitz continuous, Equation~\eqref{eq:sdegamma} even has a unique strong solution, see~\cite[Theorem~V.6]{Protter2004}.

We are ultimately interested in asymptotic properties of our Bayesian procedure, detailed below, for estimating the volatility function $\sigma$. To that end, we need a sufficiently rich set of observations and this can be accomplished by the scaled observation process in \eqref{eq:sden} below, instead of $X$ satisfying \eqref{eq:sdegamma}. So we consider the process $X^n$ given as the solution to
\begin{equation}\label{eq:sden}
\dd X^n_t=\frac{1}{n}\sigma(X^n_{t-})\,\dd L_t,\,\, X^n_0=0.  
\end{equation}
The scaling factor $\frac{1}{n}$ causes, for large values of $n,$ a `slow growth' of the process $X^n$ and `long times' to stay around certain levels.
We will thus assume that $X^n$ is observed on a long time interval $[0,T^n]$, where $T^n\to\infty$. In Section~\ref{section:asymptotics} we will specify $T^n$ and derive asymptotic results for $n\to\infty$.

\subsection{Likelihood}

For a Bayesian analysis we need to convert a prior distribution into a posterior one, and for this a likelihood ratio is needed. Let us first sketch the set-up. 

Let $\fpspace$ be a filtered probability space and let $(L_t)_{t\geq 0}$ be a gamma process adapted to $\mathbb{F}$, whose L\'evy measure admits the density $v$ given by \eqref{eq:v}.
Assume that $X$ is a (weak) solution to \eqref{eq:sdegamma},
and assume that $X$ is observed on an interval $[0,T]$. We denote by $\pp^\sigma_T$ its law, a probability measure on the Skorohod space $D[0,T]$ with the Skorohod $J_1$-topology and the $\sigma$-algebra generated by it. In agreement with this notation we let $\pp^1_T$ be the law of $X$ when $\sigma\equiv1$, in which case $X_t=L_t,\, t\in [0,T]$. The measure $\pp^1_T$ will serve as a reference measure. Recall \eqref{eq:v} and let
\[
Y(t,x):=\frac{1}{\sigma(X_{t-})}v\Bigl(\frac{x}{\sigma(X_{t-})}\Bigr)/v(x).
\]
In \cite{bgssweak}  the following proposition has been proven.

\begin{proposition}\label{prop:z}
Let $\sigma$ be a positive locally bounded measurable function on $[0,
\infty)$ such that \eqref{eq:sdegamma} admits a weak solution that is unique in law. Assume moreover that $\sigma$ is lower bounded by a strictly positive constant. Let $T$ be a  finite (stopping)  time $T>0$. Then the laws $\pp^\sigma_T$ and $\pp^1_T$ are equivalent on $\cf^X_T$ and the corresponding Radon-Nikodym derivative $Z_T$  
has the explicit representation: 
\begin{equation}\label{eq:explicit}
Z_T=\exp\left(\int_{0}^{T}\int_{0}^{\infty}\log Y(t,x)\,\mu^X(\dd x,\dd t)-  \int_0^T\int_{0}^{\infty}(Y(t,x)-1)v(x)\,\dd x\,\dd t  \right),
\end{equation}
where $\mu^X$ is the jump measure of $X$ and both double integrals are a.s.\ finite. 
\end{proposition}
In the next section we will present, under additional assumptions, an explicit expression for $Z_T$.

\subsection{Prior and posterior}\label{section:pp}

We  model $\sigma$ as a piecewise constant function:
\begin{eqnarray}
\label{eq:sigma_pc}
\sigma(x)=\sum_{k=1}^K \xi_k \mathbf{1}_{B_k}(x)
\end{eqnarray}
for bins $B_1=[0,b_1]$,
$
B_k=(b_{k-1},b_k], k=2,\ldots,K-1,
$
and  $B_K=(b_{K-1},b_K)$, with appropriately chosen increasing sequence of bin endpoints $\{b_k\}$ and the bin number $K$. 
The $\{\xi_k\}$ are positive numbers (later on positive random variables). Although we use \eqref{eq:sigma_pc} for our model, we emphasize that the `true' $\sigma$ does not need to be piecewise constant. As a final remark we note that when $\sigma$ is given by \eqref{eq:sigma_pc}, \eqref{eq:sden} still has a unique solution, obtained as a concatenation of stopped gamma processes. Under the specification~\eqref{eq:sigma_pc}, the general likelihood ratio $Z_T$ of \eqref{eq:explicit} gets an explicit form.

\begin{proposition}\label{cor:lik}
Suppose that $X^n$ is given by \eqref{eq:sden} with \(\sigma\) given by \eqref{eq:sigma_pc}. Let \(\tau^n_{k}=\inf\{t\geq 0: \, X^n_t\geq b_k\}\), \(k=1,\ldots,K,\)
and write $T=\tau^n_K$. Then
\begin{eqnarray*}
\frac{\dd\pp^\sigma_T}{\dd\pp^1_T}=\exp\left\{ \beta\sum_{k=1}^{K}\left[1-n\xi_{k}^{-1}\right](X^n_{\tau^n_{k}}- X^n_{\tau^n_{k-1}})-\alpha\sum_{k=1}^{K}\left(\tau^n_{k}-\tau^n_{k-1}\right)\log(\xi_{k}/n)\right\}. 
\end{eqnarray*}
\end{proposition}
To assign a prior distribution, we take the $\xi_k$ to be independent random variables with an inverse gamma distributions for each of them, that is \(\xi_{k}\sim \ig(\alpha_k,\beta_k)\), and thus have marginal densities $p_k(x)=\beta_k^{\alpha_k}\Gamma(\alpha_k)^{-1}x^{-\alpha_k-1}\exp(-\beta_k/x)$. 
Therefore, we possibly have to extend the original probability space to carry the $\xi_k$ as well, taking into account that $L$ and the $\xi_k$ have to be independent. We can then use Proposition~\ref{cor:lik}, where $\pp^\sigma_T$ is to be interpreted as the conditional law of $X^n$ on $[0,T]$ given the $\xi_k$. It turns out that the posterior distributions of the $\xi_k$ are again of inverse gamma type. 

\begin{lemma}
Let $x_1,\ldots,x_K$ denote realisations of $\xi_1,\ldots\xi_K$ and let $X^n$ stand for the path $X_t^n,$  $t\in [0,T]$. With $T=\tau^n_K$, we then have from Proposition~\ref{cor:lik} that the posterior joint density of $(\xi_1,\ldots,\xi_K)$ is given by
\[
\Pi(x_1,\ldots,x_K\mid X^n)\propto\prod_{k=1}^{K}\exp\left(-(n\beta (X^n_{\tau^n_{k}}- X^n_{\tau^n_{k-1}})+\beta_{k})x_{k}^{-1}\right)x_{k}^{-\alpha(\tau^n_{k}-\tau^n_{k-1})-\alpha_{k}-1}.
\]
It follows that the $\xi_k$ are independent under the posterior distribution, and
\begin{equation}\label{eq:igposterior}
\xi_k \mid X^n \sim \operatorname{IG}(\alpha \Delta \tau^n_k + \alpha_k,  n\beta \Delta X^n_{\tau^n_k} + \beta_k),
\end{equation}
where $\Delta \tau^n_k=\tau^n_{k}-\tau^n_{k-1}$, $\Delta X^n_{\tau^n_k}=X^n_{\tau^n_{k}}- X^n_{\tau^n_{k-1}}$.
\end{lemma}

\subsection{Asymptotics}\label{section:asymptotics}

We consider the process $X^n$ satisfying \eqref{eq:sden}, but the standing assumption from here on is that $\sigma:[0,\infty)\to[0,\infty)$ is H\"older continuous, that is, there are constants $H\geq 0$ and $0<\lambda\leq 1$ such that for all $x,y>0$ it holds that $|\sigma(x)-\sigma(y)|\leq H|x-y|^\lambda$.  Moreover, $\sigma$ is assumed to be bounded from below by a positive constant $\underline{\sigma}$.  
Here are some further assumptions. 
\begin{itemize}
\item
The number of bins and their width depend on $n$. So we write $B^n_k=(b^n_{k-1},b^n_k]$, $k=1,\ldots,K$, $K=K_n$. We assume equidistant bins. Let $b_K$ be the endpoint of the last bin, assumed to be fixed. We take the other bin boundaries $b^n_k$ as $b^n_k=\frac{b_Kk}{K}$, $k=1,\ldots,K.$ A given $x\in (0,b_K]$ then belongs to bin $B^n_k$ with $k=k_n(x)=\lceil \frac{Kx}{b_K}\rceil$.
\item
We assume the number of bins $K=K_n\asymp n^\kappa$ for $0<\kappa<1$. Then, given also the above assumption on the $b^n_k$, one has $\Delta b^n_k:=b^n_k-b^n_{k-1}\asymp n^{-\kappa}$ for all $k$ and $x\in B_k$. Note that the bin widths decrease to zero and that the number of bins increases to infinity, with a rate depending on the parameter $\kappa$ (to be specified below). This behaviour is needed to have better and better approximations of the true function $\sigma$ with piecewise constant functions.
\item
We observe the process $X^n$ until it crosses the last bin, i.e.\ until time $T^n=\tau^n_K$. In \cite{bgsslevy} it has been shown that the random time $T^n$ is roughly proportional to $n$.
\end{itemize}
Although $\sigma$ is continuous, we assume it to be apriori  piecewise constant. That is, parallel to \eqref{eq:sigma_pc} but note that the bins now depend on $n$, as
\begin{equation}\label{modelsigma}
\xi^n(x)=\sum_{k=0}^K\xi_k\one_{B^n_k}(x),
\end{equation} 
where the $\xi_k$ are assigned the inverse gamma prior distributions as in Section~\ref{section:pp}. Consequently, the posterior distributions of the $\xi_k$ are as in \eqref{eq:igposterior}.

Here is the main result, it says that the contraction rate of the posterior distribution is $n^{-\lambda/(2\lambda+1)}$, $\lambda$ being the H\"older exponent of $\sigma$. We use the notation $\Pi_n$ to denote posterior probabilities given the observations of $X^n$ on the time interval $[0,T^n]$.

\begin{theorem}
\label{thm:rates}
Assume the model with  volatility  \eqref{modelsigma} whereas the true volatility function $x\mapsto \sigma(x)$ is H\"older continuous of order $\lambda\leq 1$ and bounded from below. Let the bin sizes shrink proportional to $n^{-\frac{1}{2\lambda+1}}$ (i.e.\ $\kappa=\frac{1}{2\lambda+1}$), and let $(m_n)$ be any sequence of real numbers diverging to infinity. Then, for $n\to\infty$ 
\[
\sup_{x\in [0,b_K]}\Pi_n\Bigl(|\xi^n(x)-\sigma(x)|>m_n n^{-\frac{\lambda}{2\lambda+1}}\Bigr)\to 0 \mbox{ in $\pp^\sigma_{T^n}$-probability}.
\]
Moreover, this contraction rate is minimax optimal.
\end{theorem}

\subsection{Numerical illustration}

\begin{figure}
\begin{center}
\includegraphics[width=0.7\linewidth]{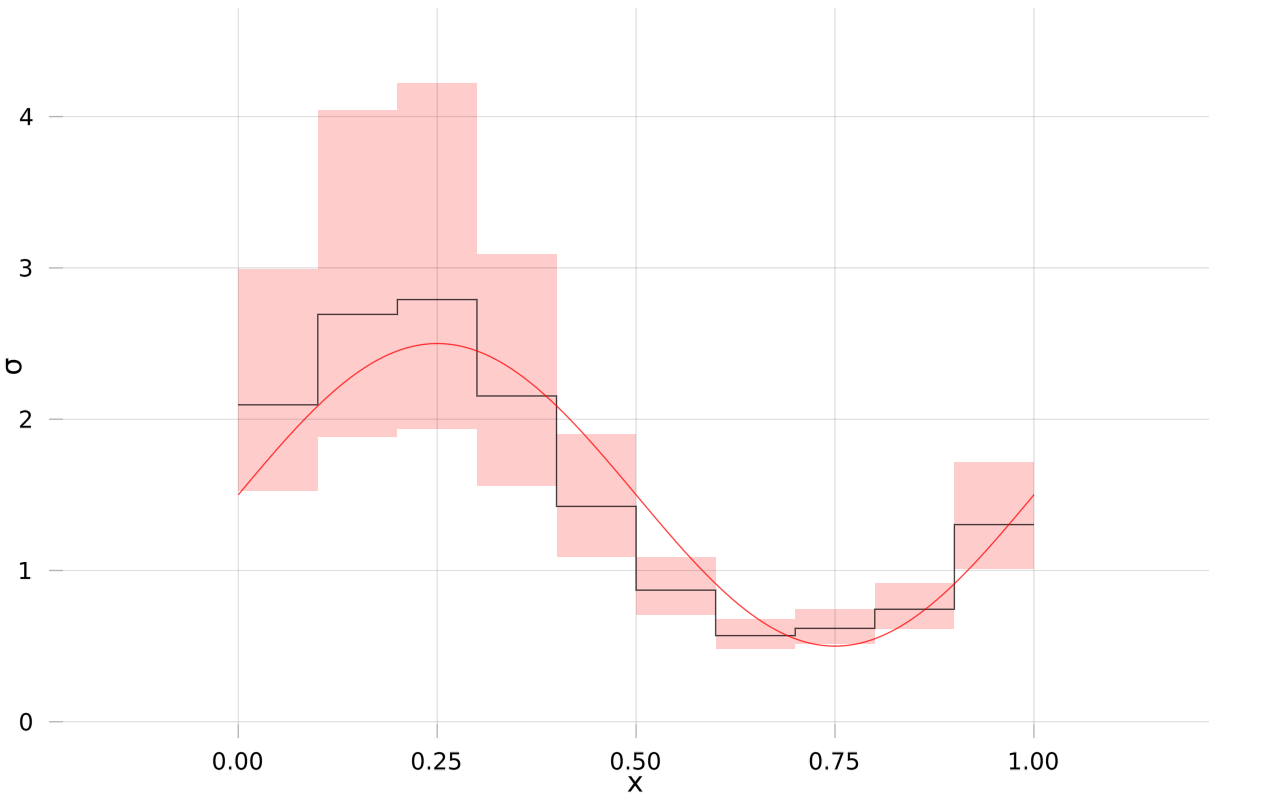}
\end{center}
\caption{Red: true volatility $\sigma_0$ as a function of $x$. Shaded: marginal $90\,\%$-posterior credible band for the piecewise constant posterior $\sigma$. Black: marginal posterior median. }
\label{fig:ex1post}
\end{figure}

Here we graphically illustrate our inferential results. In agreement with the notation pertaining to the asymptotic regime detailed in \eqref{eq:sden}
we took $X^n$ as the solution to the L\'evy SDE  \eqref{eq:sdegamma} with the volatility function
\begin{eqnarray}
\label{eq:sigma0}
\sigma_0(x) = \frac{3}{2}+\sin(2\pi x)
\end{eqnarray}
and scaled it with $1/n$, $n=500$. 
We have partitioned the unit interval into $10$ bins setting $b_k = k \delta x$, $\delta x = 0.1,$ $k=1,\ldots,10,$ and used the piecewise constant prior of the form $\sigma(x)=\sum_k \xi_k \mathbf{1}_{B_k}(x),$ where $\xi_1,\ldots,\xi_{10}$ are i.i.d.\ random variables with the inverse gamma $\operatorname{IG}(\alpha_k,\beta_k)$ distribution with $\alpha_k=\beta_k=0.1$. The posterior is given in closed form  by \eqref{eq:igposterior}. Figure~\ref{fig:ex1post} contrasts the corresponding marginal  $90\,\%$-posterior credible bands for $\sigma$ with the true volatility function $\sigma_0$ from 
\eqref{eq:sigma0}.

\subsection{Closing remarks}

There are many extensions conceivable. The inclusion of a drift in Equation~\eqref{eq:sdegamma} is one possibility, another one is to move from a continuous record of observations to a discrete time setting. See \cite{bgsslevy} for short discussions of such possible extensions and also for an extensive treatment of an application to real data obtained from the North Greenland Ice Core Project. The present section concentrated on the theoretical side.

\section{Gamma-type L\'evy subordinators}\label{section:gamma}

Given discrete time observations over a growing time interval, we consider in this section a nonparametric Bayesian approach to estimation of the L\'evy density of a L\'evy process belonging to a flexible class of infinite activity subordinators. We establish that our nonparametric Bayesian procedure is consistent: in the low frequency data setting, with equispaced in time observations and intervals between successive observations remaining fixed, the posterior asymptotically,  as the sample size $n\rightarrow\infty$, concentrates around the L\'evy density under which the data have been generated. Finally, we test our method on a classical insurance dataset. This section is based on \cite{belomestny18}.

\subsection{Problem formulation}

Consider a univariate L\'evy process $X=(X_t\colon t\geq 0)$ with L\'evy triplet $(\gamma,A,\nu).$ By the L\'evy-Khintchine formula, see Theorem 8.1 in \cite{Sato2013}, the characteristic function $\phi_{X_1}$ of $X_1$ admits the unique representation of the type
\[
\phi_{X_1}(z)=\exp\left( \mathrm{i} \gamma z-\frac{1}{2} A z^2+\int_{\mathbb{R}} \left( e^{\mathrm{i}zx}-1-\mathrm{i}zx\one_{|x|\leq 1} \right) \nu (\dd x) \right),
\]
and the triplet $(\gamma,A,\nu)$ uniquely determines the law of the process $X.$ We assume that the L\'evy measure \(\nu\) admits the representation
\begin{equation}\label{levy}
\nu(\dd x)=\frac{\beta}{x}e^{-\alpha x-\theta(x)} \dd x, \quad x > 0,
\end{equation}
where $\alpha$ and $\theta\colon[0,\infty)\rightarrow \mathbb{R}$ are parameters to be estimated, while $\beta$ is a known or unknown parameter. 
This entails that $\theta$ and $\alpha$ are such that $\nu([1, \infty))$ is finite. We also assume that $A = 0$ and
\begin{equation}
\label{cond.gamma}
\gamma=\int_{0\leq x \leq 1} x\nu(\dd x) <\infty.
\end{equation}
It follows that  $X$ is a pure jump process with non-decreasing sample paths,
or put another way, a subordinator with zero drift, cf.\ Sections 2.6.1--2.6.2 in \cite{kyprianou14}. As in our case the L\'evy triplet $(\gamma,A,\nu)$ is entirely determined by $\nu$, we will write $\pp_\nu$ for the law of $X$.
One may call this class of L\'evy processes Gamma-type subordinators,
because $X$ is a Gamma process when $\theta \equiv 0$, but we prefer to simply refer to it as $\theta$-subordinators.

Assume that the process $X$ is observed at discrete time instances $0 = t_0 <t_1 <\dots < t_n = T,$ so our observations are $X^{(n)} = (X_{t_i} \colon  i \in \{0,\ldots, n\})$. Our aim is nonparametric Bayesian estimation for the parameter triple $(\alpha,\beta,\theta)$. This requires specification of the likelihood and the prior in our model, that are next combined via Bayes' formula to form the posterior distribution. This latter encodes all the necessary inferential information within the Bayesian setup.

\subsection{Likelihood}

We consider the problem with $\beta$ known and fixed in the present treatment. This choice has been made because of two reasons. First, the parameter $\beta$ plays a role similar to the diffusion coefficient $\sigma$ in a stochastic differential equation driven by a Wiener process; different values of $\sigma$ there result in mutually singular laws of a continuous record of observations, and similar singularity is valid for $\theta$-subordinators. Derivation of nonparametric Bayesian asymptotics for the stochastic differential equations (all of which is a recent work) historically proceeded from the assumption of a known $\sigma$ to the one where $\sigma$ is unknown and has to be estimated. In that sense, the fact that we assume $\beta$ is known at this stage does not appear unexpected or unnatural.
There is also a technical reason to assume $\beta$ to be known.
It assists in deriving useful bounds on the Kullback-Leibler and Hellinger distances between marginals of $\theta$-subordinators under different L\'evy triplets, which in general, is the key to establishing consistency properties of nonparametric Bayesian procedures. 
For a treatment of the statistical problem with observations in discrete time when also $\theta$ is unknown, we refer to \cite[Section~6]{belomestny18}, and for now, we maintain the assumption of known $\beta$.

All processes and their laws in this section are restricted to the time interval $[0, T]$ for a fixed $T > 0$.
Note that for any two L\'evy measures $\nu$ and $\nu_0$ given by \eqref{levy} with parameters $\beta,\alpha, \theta$ 
and $\beta,\alpha_0,\theta_0,$ respectively, provided $\theta(0) = \theta_0(0) = 0$ and both functions $\theta$ and $\theta_0$ are Lipschitz continuous in some neighbourhood of zero, we have 
$\nu$ and $\nu_0$ are equivalent. Moreover,
\begin{equation}
\label{cond.nu}
d^2_\cH(\nu, \nu_0)=\frac{1}{2}\int\limits_{(0,\infty)} ( \sqrt{\dd \nu} - \sqrt{\dd \nu_0} )^2<\infty,
\end{equation}
where $d_\cH(\cdot,\cdot)$ is the Hellinger distance between two (infinite) measures. By assumption \eqref{cond.gamma} and property \eqref{cond.nu}, together with Theorem~33.1 in \cite{Sato2013}, it follows that the laws $\pp_{\nu}$ and $\pp_{\nu_0}$ of $X =(X_t\colon t\in [0,T])$ are equivalent. Furthermore, Theorem~33.2 in \cite{Sato2013} implies that a.s.\
\begin{equation*}
U_T=\log\left( \frac{\dd \pp_{\nu}}{\dd \pp_{\nu_0}}\big(X\big)\right) = \sum_{(s,\Delta X_s)\in (0,T]\times\{ \Delta X_s>0 \}} \phi(\Delta X_s) - T \int\limits_{\mathclap{(0,\infty)}} (e^{\phi(x)}-1)\nu_0(\dd x) ,
\end{equation*}
where $\Delta X_s=X_s-X_{s-},$ and
\[
\phi(x)=\log \left(\frac{\dd \nu }{\dd \nu_0}(x)\right) = -(\alpha x + \theta(x) - \alpha_0 x - \theta_0(x)),  \quad x>0.
\]
We can also write the log-likelihood ratio $U_T$ as
\[
U_T=\int_{(0,T]} \int_{{(0,\infty)}} \phi(x) \mu(\dd s,\dd x) - T\int\limits_{\mathclap{(0,\infty)}} (\nu-\nu_0)(\dd x),
\]
where the jump measure $\mu$ is defined by
\[
\mu((0,t]\times B)= \#\left\{ s\colon (s,\Delta X_s)\in (0,t]\times B\right\}
\]
for any Borel subset $B$ of $(0,\infty)$.
We can view  $\pp_{\nu_0}$ as the dominating measure for $\pp_{\nu}$. From an inferential point of view the specific choice of the dominating measure is immaterial.  A convenient choice of $\nu_0$ for the theoretical development in Section~\ref{sec:asymp} is to actually take $\nu_0$ to be the `true' L\'evy measure $\nu_0$ with parameters $\alpha_0$ and $\theta_0$.

\subsection{Prior}

To define the prior, we consider a subclass of processes defined in  \eqref{levy}, where
the parameter $\theta$ in the L\'evy measure $\nu$ has the following form. Fix a sequence
$0 <b_1<\cdots<b_{N}<\infty$,
set for convenience  $b_{0} = 0$ and $b_{N+1}=\infty,$ and define bins $B_k$ by
\[
B_k=[b_{k},b_{k+1}), \quad k=0,\ldots, N.
\]
Given bins $B_k,$ assume the function $\theta$ is, unlike those in previous sections, piecewise linear, i.e.,
\begin{equation}
\label{piecewise}
\theta(x)=\sum_{k=1}^N (\rho_k+\theta_k x) \one_{B_k},
\end{equation}
where $\rho_k\in\mathbb{R},$ $k=1,\ldots,N,$ $\theta_k\in\mathbb{R},$ $k=1,\ldots,N,$ and $\theta_N>-\alpha.$ Together with $\alpha,$ the parameter $\theta_k$ determines the slope of the function $\theta(x) + \alpha x$ on the bin $B_k,$ while $\rho_k$ gives the intercept. The process $X$ with the law $\pp_{\nu}$ can be viewed as a Gamma process with rate parameter $\alpha$ and shape parameter $\beta$,  subjected to local deviations in the behaviour
of jumps of sizes falling in bins $B_k$ compared to that of a Gamma process. The parameters $\theta_k,\rho_k$ quantify the extent of these local deviations on the bin $B_k$.
\par
We model $\alpha,\theta_k,\rho_k$ as independent under the prior distribution. Note such a prior on $\alpha,\theta_k,\rho_k$ implicitly defines a prior on the L\'evy measure $\nu$ as well.
The specific form of the prior is not crucial for many arguments that follow, but  is convenient computationally. In fact, theoretical results in Section~\ref{sec:asymp} can be derived for other series priors as well. However, the local linear structure in \eqref{piecewise} (which also means that  the prior could be rewritten as a series prior with basis functions having compact support) is important to derive some simple update formulae below.
For a realisation $\nu$ from the implicit prior on $\nu$ as above in the present section, let us work out the integral
\[
\nu(B_k)=\int_{b_{k}}^{b_{k+1}} \frac{\beta}{x} e^{-(\alpha + \theta_k) x-\rho_k}\dd x,
\]
which enters the expression for the likelihood. To that end, remember the definition of the exponential integral,
$
E_1(z)=\int_z^{\infty} t^{-1}{e^{-t}}\dd t,
$
see, e.g., \S 15.09 in \cite{jeffreys1999methods} for its basic properties. Then a change of the integration variable gives
\begin{equation}
\label{nu.bk}
\nu(B_k)=\beta e^{-\rho_k} \{ E_1((\theta_k + \alpha) b_{k}) - E_1((\theta_k + \alpha) b_{k+1}) \}, \quad k=1,\ldots,N.
\end{equation}
Observe that $\nu(B_N)=\beta e^{-\rho_N} E_1((\theta_k + \alpha) b_{N}).$ Similar to the case of $\nu$,
\[
\nu_0(B_k)=\beta\{ E_1(\alpha_0 \,b_{k}) - E_1(\alpha_0 \,b_{k+1}) \}, \quad k=1,\ldots,N.
\]
Also here remark that $\nu_0(B_N)=\beta E_1(\alpha_0 \,b_{N}).$ For future reference, note that for any $\alpha,\alpha^{\prime},$
\begin{equation}
\label{frullani}
\lim_{x \to 0} \{ E_1(\alpha x) -  E_1(\alpha' x) \}= \log\left( \frac{\alpha'}{\alpha }\right),
\end{equation}
which follows from the formula for Frullani's integral, see \S 12.16 in \cite{jeffreys1999methods}.

\subsection{Asymptotics}
\label{sec:asymp}
Let the process $X$ be observed at equidistant times $t_i$, $i=1,\ldots,n$. Without loss of generality we assume that the observations are $X_1,\ldots,X_n.$ We pass to the increments $Z_i=X_i-X_{i-1},$ $i=1,\ldots,n,$ and denote our observations by $\mathcal{Z}_n=\{Z_1,\ldots,Z_n\}.$ We assume that under the true L\'evy density $v_0,$ $Z_1 \sim \mathbb{Q}_{v_0}$. In general, $\mathbb{Q}_{v}$ will stand for the law of the increment $Z_1$ under the L\'evy density $v.$ Furthermore, we introduce the law  $(X_t\colon t\in[0,1]) \sim \mathbb{R}_{v_0}$ under the true L\'evy density $v_0.$ The law of this path under the L\'evy density $v$ will be denoted by $\mathbb{R}_v.$ For  asymptotic results, we will let the number of bins $N$ depend on the sample size $n$, and write $N_n$ instead. We postulate a condition on the underlying model and on the chosen prior distribution.

\begin{cnd}
\label{cnd:truth}
\emph{
Let the function $\theta_0$  have a compact support on the interval $[\underline{b},\overline{b}]$ where the boundary points $0<\underline{b}<\overline{b}<\infty$ are known, $\|\theta_0\|_\infty<\bar\theta$, and suppose  $\theta_0$ is $\lambda$-H\"older continuous, that is, $|\theta_0(x)-\theta_0(y)|\leq L|x-y|^{\lambda}$ ($\lambda\in (0,1]$, $L>0$). Suppose also that $\alpha_0\in [\underline{\alpha},\overline{\alpha}]$ with known boundary points $0<\underline{\alpha}<\overline{\alpha}<\infty$. Finally, assume that the parameter $\beta_0$ is known (and, without loss of generality, equal to $1$.
}\end{cnd}

\begin{cnd}
\label{cnd:prior}
\emph{
The prior $\Pi_n$ is defined on a special class of L\'evy densities, $V_n$. These are the densities that on the bins $B_k=(b_{k-1},b_k]$, $k=1,\ldots,N$, $b_0=0$, $b_1=\underline{b}$, $b_N=\overline{b}$, have the form $v(x)=\frac{\beta_0}{x}\exp(-\alpha x -\theta_k(x))$, with $\theta_k(x)=\rho_k+\theta_kx$, with the special choice $\rho_0=\theta_0=0$ and $\beta_0=1$. So, with the above notation,
\[
V_n=\left\{v\colon v_{|B_k}(x)=\frac{\beta_0}{x}\exp(-\alpha x -\theta_k(x)), \, k=1,\ldots,N\right\}.
\] 
Under the prior distribution, the coefficients $\theta_i,$ $i=1,\ldots N-1,$ are independent and uniformly distributed on the known interval $[-\overline{\theta},\overline{\theta}]$, $\overline\theta>0$.  Likewise, the coefficients $\rho_i,$ $i=1,\ldots,N-1,$   are independent uniform on the interval $[-\overline{\theta},\overline{\theta}],$ whereas $\alpha$ is uniform on $[\underline{\alpha},\overline{\alpha}]$, $\overline\alpha>0$. We assume a priori independence of all $\rho_i$, $\theta_i$ and $\alpha$. Implicitly, this defines a prior $\Pi_n$ on the class of L\'evy densities $V_n$, which are realisations from the prior.
}
\end{cnd}
Theorem~\ref{thm:consistency} below implies that our Bayesian procedure is consistent in probability; this in turn implies the existence of consistent Bayesian point estimates, see, e.g., \cite{ghosal2000convergence}, pp.~506--507.
We use the notation $\Pi_n(\dd v\mid \mathcal{Z}_n)$ for the posterior measure. Also, $\mathbb{Q}_{v_0}^n$ denotes the law of the sample $\mathcal{Z}_n$ under the true L\'evy density $v_0$. 
\begin{theorem}
\label{thm:consistency}
Assume that Conditions~\ref{cnd:truth} and~\ref{cnd:prior} hold and that $N_n\to\infty$ and $N_n/n\to 0$ as $n\to\infty$.
Let $d_{\mathcal{H}}$ be the Hellinger metric. Then, for any fixed $\epsilon,\varepsilon'>0,$
\[
\mathbb{Q}_{v_0}^n \left( \Pi_n(v\colon d_{\mathcal{H}}(\mathbb{Q}_{v_0},\mathbb{Q}_v) > \epsilon \mid \mathcal{Z}_n) > \varepsilon' \right) \rightarrow 0
\]
as $n\rightarrow\infty.$ 
\end{theorem}

\subsection{Numerical results}

Over the last two decades there has been an increasing interest in applying Bayesian methods to insurance problems, see, e.g., \cite{hong2017review} and references therein. \cite{hong2018dirichlet} apply a Dirichlet process mixture prior to model the density of insurance claim sizes, and provide motivation for using a nonparametric Bayesian approach in the actuarial science. In this section we will apply our Bayesian approach to the Danish data on large fire losses.  This dataset is a standard test example in extreme value theory, and from that point of view it has been a subject of several deep studies, such as \cite{mcneil1997estimating} and \cite{resnick1997discussion}. Our goals here are more modest, and aim at demonstrating the facts that firstly, $\theta$-subordinators can be potentially used to capture some aggregate features of the Danish data on large fire losses, and secondly, statistical inference for real data modelled through such processes can be successfully performed using the Bayesian methodology developed here. 
\subsubsection*{Data description and visualisation}

\begin{figure}
\includegraphics[width=0.45\textwidth]{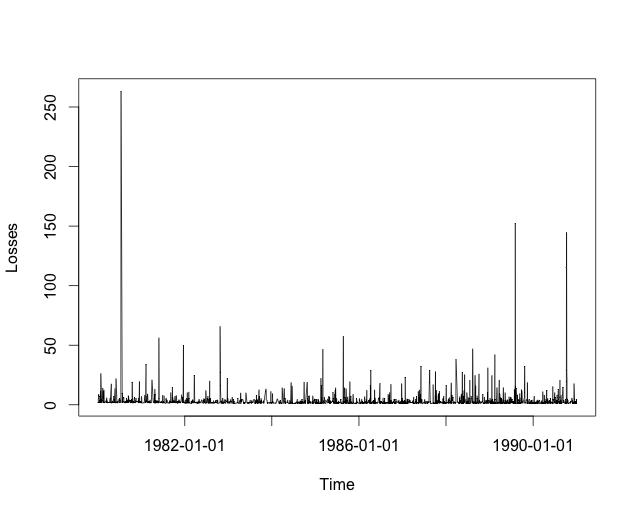}
\includegraphics[width=0.45\textwidth]{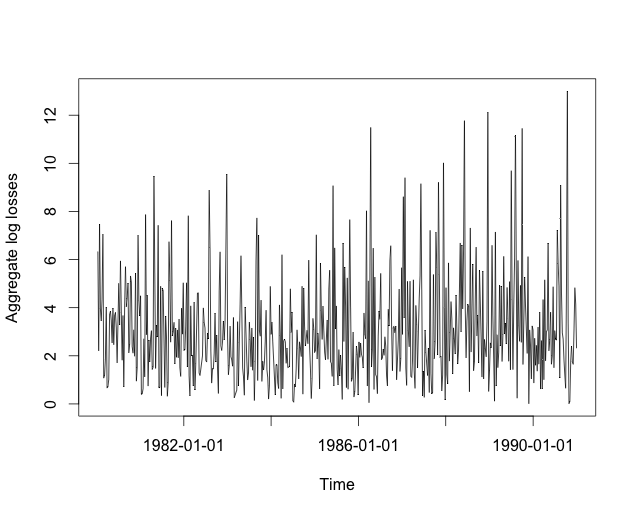}
\caption{Danish data on large fire losses. {\em Left}: original daily data (the unit is one million DKK). {\em Right}: logarithmically transformed and aggregated data.}
\label{fig:danishdata}
\end{figure}

A succinct description of the Danish data on large fire losses can be found on p.~298 in \cite{embrechts2013modelling}. The dataset (scaled for privacy reasons) comprises 2167 fire losses (adjusted suitably for inflation to reflect the 1985 values) in Denmark over the 10 year period starting on 6 January 1980 and ending on 30 December 1990, that exceed in size one million DKK, and that were registered by Copenhagen Reinsurance. The rationale for thresholding losses at one million DKK is given in \cite{mcneil1997estimating}, pp.~119--120, and consists in the fact that in practice it is virtually impossible to collect exhaustive data on small losses: insurance is typically provided against significant losses, while small losses are dealt with by insured parties directly.
The data can be accessed through the {\bf QRM} package in {\bf R} under the name {\texttt{danish}}. The time plot of the data is given in the left panel of Figure~\ref{fig:danishdata}.

According to the exploratory analysis of the transformed data, the data can be modelled as an i.i.d.\ sequence that follows a Gamma-like distribution, but perhaps is not genuinely Gamma. This suggests a possibility of using a $\theta$-subordinator to model the data.

\subsubsection*{Modelling fire losses with a $\theta$-subordinator}

\begin{figure}[htbp]
\begin{center}
\includegraphics[width=\textwidth]{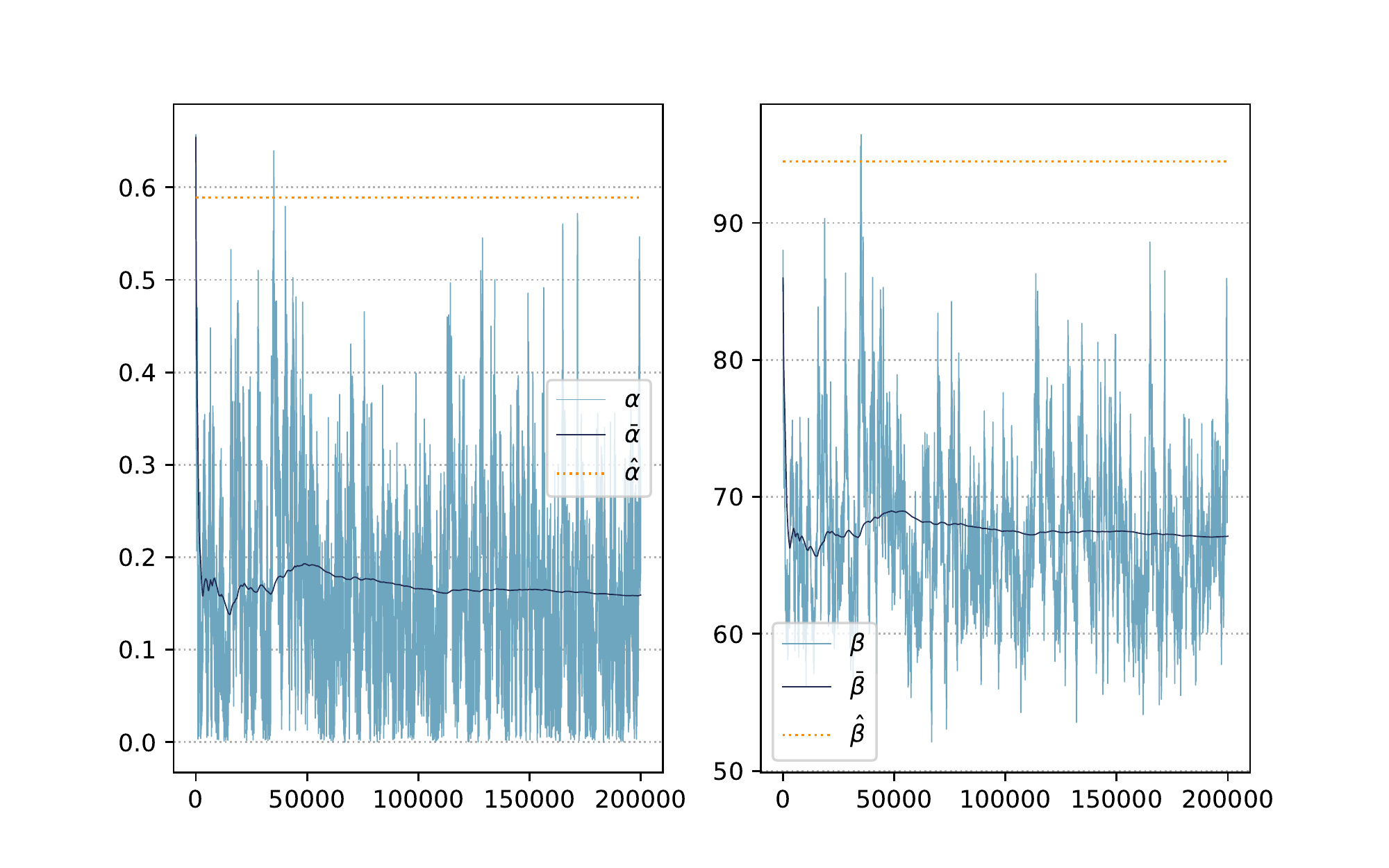}
\caption{Trace plots of the parameters $\alpha$ and $\beta$ for the fire loss data. Left: trace and running average of samples of $\alpha$.
(The latter indicated by $\bar \alpha$.) 
The maximum likelihood estimate $\hat\alpha$ of $\alpha$ using a Gamma process model is marked as the dotted yellow line.
Right: trace and running average of samples of $\beta$. (The latter indicated by $\bar \beta$.) The maximum likelihood estimate $\hat\beta$ of $\beta$ using a Gamma process model is marked as the yellow dotted line.
}
\label{danish:traceplot1b}
\end{center}
\end{figure}

\begin{figure}[htbp]
\begin{center}
\includegraphics[width=\textwidth]{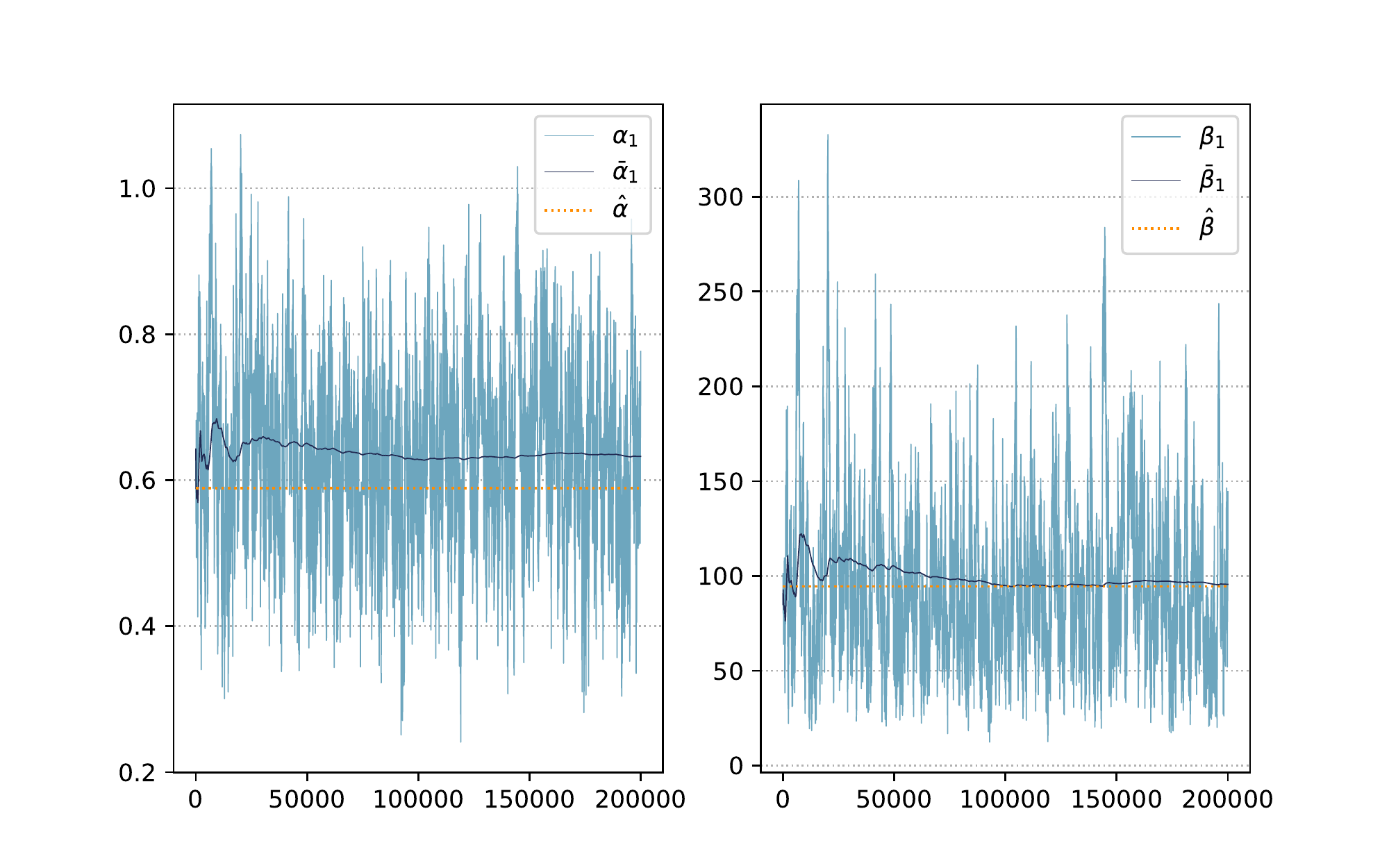}
\caption{Trace plots of the parameters used for the bin $(b_1, \infty)$ for the fire loss data. Left: trace and running average of the samples of $\alpha_1$.
The maximum likelihood estimate of $\alpha$ using a Gamma process model is marked as yellow line.
Right: trace and running average of the samples of $\beta_1$. 
The maximum likelihood estimate of $\beta$ using a Gamma process model is marked as yellow line.
}
\label{danish:traceplot2b}
\end{center}
\end{figure}

\begin{figure}[htbp]
\begin{center}
\includegraphics[width=\textwidth]{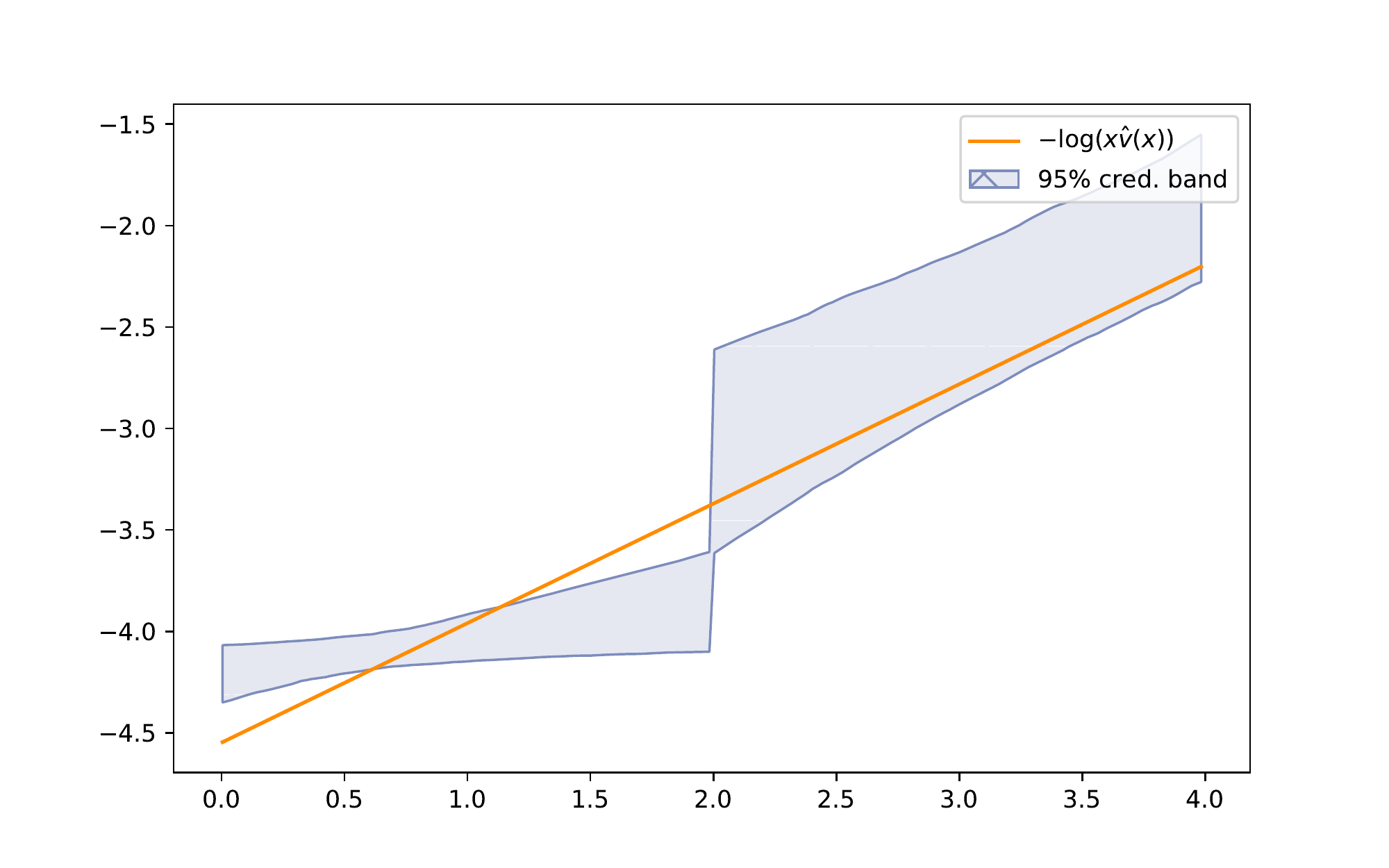}
\caption{Marginal Bayesian credible bands for the fire loss data for the function $-\log(x v(x))$ based on all samples.
Yellow: maximum likelihood estimate $-\log(x \hat v(x))$ assuming a Gamma process.
}
\label{danish:bandsb}
\end{center}
\end{figure}

Because the sample size is much smaller compared to our simulation examples, we chose $N=1$ corresponding to a single grid point $b_1 = 2$ and four parameters $\alpha$, $\beta$, $\theta_1$, $\rho_1$. In  order to improve mixing of the chain, we use a reparametrization
$\alpha_1 = \alpha + \theta_1$, $\beta_1 = \beta\exp(-\rho)$, and work with four parameters  $\alpha$, $\beta$, $\alpha_1$, $\beta_1$, so that 
\[
v(x) = \begin{cases}\frac{\beta}x\exp(-\alpha x) &x \le b_1,\\
\frac{\beta_1}x\exp(-\alpha_1 x) &x > b_1.
\end{cases}
\]
A priori we equip these four parameters with Gamma distributions, with mean 0.75 and variance $0.36$ for  the parameters $\alpha, \alpha_1$, and mean $90$ and variance $2500$ for the parameters $\beta$, $\beta_1$, and in such a way that they become independent random variables.
In the data augmentation step we take intermediate points at distance $0.0192$, corresponding to $m = 1000$. 

For the parameter updates we took independent Gaussian innovations with standard deviations
$\sigma_\alpha = \sigma_{\alpha_1} = 0.03$, $\sigma_\beta = 1$ and $\sigma_{\beta_1} = 6$, respectively.
In the Gibbs sampler in each step new Gamma bridges are proposed  in the data augmentation step, followed by a parameter update Metropolis-Hastings step cycling through updates of $\beta$ in the first and second and the other parameters jointly in each of the remaining three of in total 5 stages. 
With these choices, the chains mix sufficiently well. 
The MCMC algorithm was ran for $200\,000$ iterations.  
Figure~\ref{danish:traceplot1b}
shows trace plots and running averages of the posterior samples of the 
parameters $\alpha$ and $\beta$, whereas Figure~\ref{danish:traceplot2b} shows similar plots for the parameters $\alpha_1$ and $\beta_1$.

Figure~\ref{danish:bandsb} shows the $95\%$ marginal Bayesian credible band for the function $\theta(x)+\alpha x$ contrasted to the maximum likelihood estimate that assumes the observations come from a Gamma process. This plot suggests that modelling the losses with a Gamma process leads to overestimation of the number of small jumps and possibly of large jumps too; however, more data is necessary to make a definitive statement (unfortunately, as observed in \cite{chavez2016extreme}, it is difficult for academia to gain access to the insurance data).
In connection to this, we note that a difference in the estimates of the rate of decay of the L\'evy density (value of $\alpha_1$ in the model) has serious implications of practical relevance for the assessment of the risk of very large fire losses.

\subsection{Conclusions and possible extensions}

We have considered in this section a nonparametric Bayesian approach to estimation of the L\'evy density of a L\'evy process belonging to a flexible class of infinite activity subordinators and  established that this procedure is consistent. We have also tested our method on a classical insurance dataset. 
\par
There are various directions for extensions of the model studied in this section. One may consider a class of increasing, infinite activity L\'evy processes $X$ for which the characteristic function of $X_1$ is given by
\[
\phi(z)=E[e^{\mathrm{i} zX_1}]=\int  (e^{\mathrm{i} zx}-1)\nu(dx),
\]
where $\theta(\cdot)$ is nonnegative and the L\'evy measure $\nu$ is given by
\[
\nu(\dd x) = \frac{1}{x^{1+a}}e^{-\theta(x)}\one_{(0,\infty)}(x)\dd x.
\]
This model generalizes the Gamma process with $a = 0$ and $\theta(x) \equiv\lambda$ and also covers the class of one-sided tempered stable processes.
Here one could think of estimating the stability index $\alpha$ from a Bayesian point of view.  This is difficult due to singularity of the measures induced by two L\'evy processes with different stability indices, and so this poses an interesting challenge that deserves to be explored.

\section{Summary and concluding remarks}\label{section:conclusion}

In this paper, we have reviewed some of our contributions, obtained with our collaborators, to Bayesian nonparametric function estimation. We have considered a large variety of models (volatility estimation
in diffusion models, Poisson intensity estimation, estimation for gamma-driven stochastic differential equations, estimation in the presence of microstructure noise, and L\'evy density estimation for Gamma-type subordinators), which all had a common denominator, the function to be estimated, is well approximated by piecewise constant functions (with one exception that deals with piecewise linear functions) giving rise to tractable prior distributions for Bayesian inference.

The goal of inference in all models under consideration was to obtain posterior consistency of the proposed estimators or even to quantify posterior contraction rates. 
We have seen that rougher continuous functions are more difficult to estimate than smoother ones. The degree of smoothness is in our approach measured by the H\"older exponent; the smoother the function is, the better the contraction rate of the estimation procedure is. The precise results are formulated as Theorems~\ref{thm:rate-posterior-ig}, \ref{thm:post},  and \ref{thm:rates}. 
The conclusion that emerges from these theorems is that despite the entirely different settings and the different kinds of functions (their arguments being time or space variables), the contraction rates, all depending on the H\"older exponent, are the same across all cases. A deeper analysis in terms of equivalence of statistical experiments is likely to explain the similarity of the rates. These rates were obtained under an independence assumption on the prior,  i.e., the approximating functions with constant levels on each of their bins are such that these levels are a priori independent. The setting of Section~\ref{section:gamma} was somewhat different, and we had to content ourselves with posterior consistency only, Theorem~\ref{thm:consistency}.

It also appeared that better numerical results, visible in the plots, can be obtained by relaxing the independence assumption on the prior distribution. This observation has been accomplished by using a specific Markov chain prior on the levels of the piecewise constant functions, introducing a priori dependence between the levels on adjacent bins of the piecewise constant approximating functions. 
However, the contraction rates induced by such a prior cannot improve upon the one for the prior that renders the levels independent, as the latter has been seen to be optimal.

Throughout the paper, we have given ample numerical illustrations that support the theoretical results, ranging from simulated experiments to real-world cases.


\begin{acknowledgement}
The authors gratefully acknowledge the contributions to this paper by Shota Gugushvili.
The work of D.\ Belomestny was done within the framework of the HSE University Basic Research Program.

\end{acknowledgement}

\bibliographystyle{spmpsci}

\end{document}